\documentclass{amsart}

\usepackage{amsthm,amsmath}
\usepackage{amssymb, amsfonts}
\usepackage{chngcntr}
\usepackage{color}
\usepackage{dsfont}
\usepackage{enumitem}
\usepackage{fancyhdr}
\usepackage{hyperref}
\usepackage{lipsum}
\usepackage{mathtools}
\usepackage{tikz}
\usepackage{tikz-cd}
\usepackage{ulem}
\usepackage{upgreek}
\usepackage{wasysym}

\usetikzlibrary{trees}

\newtheorem{thm}{Theorem}[section]
\newtheorem{prop}[thm]{Proposition}
\newtheorem{cor}[thm]{Corollary}
\newtheorem{lemma}[thm]{Lemma}
\newtheorem{lem}[thm]{Lemma}
\newtheorem*{thm1}{Theorem \ref{thm:Kim-dividing->coheir-dividing}}
\newtheorem*{cor2}{Corollary \ref{cor:Kim-forking >> quasi-divide}}
\newtheorem*{cor3}{Corollary \ref{cor:existence of strong coheir}}
\newtheorem*{prop4}{Proposition \ref{prop:necessary cond}}
\newtheorem*{prop5}{Proposition \ref{prop:sufficient cond}}
\theoremstyle{definition}
\newtheorem{question}[thm]{Question}
\newtheorem{fact}[thm]{Fact}
\newtheorem{rmk}[thm]{Remark}

\newtheorem{dfn}[thm]{Definition}

\newtheorem{nota}[thm]{Notation}

\newtheorem{remark}[thm]{Remark}

\newtheorem{rem/def}[thm]{Remark/Definition}
\newtheorem{notation}[thm]{Notation}

\newtheoremstyle{claim}
  {\topsep}
  {\topsep}
  {}
  {}
  {\itshape}
  {.}
  {.5em}
   {\thmname{#1}\thmnote{ (#3)}}
\theoremstyle{claim}

\newtheorem{claim}{Claim}

\newtheorem{subclaim}{Subclaim}

\newtheorem*{claim*}{Claim}

\counterwithin{figure}{thm}

\def \tree {2{^{<\omega }}}

\def \treo {\omega{^{<\omega}}}
\def \troo {\omega{^{\omega}}}
\def \trec {\lambda{^{<\kappa }}}

\def \tri {\trianglelefteq}

\def \lex {<_{lex}}
\def \len {<_{len}}
\newcommand{\trn}{%
  \mathrel{\ooalign{$\lneq$\cr\raise.21ex\hbox{$\lhd$}\cr}}}
\newcommand{\trrn}{%
  \mathrel{\ooalign{$\gneq$\cr\raise.21ex\hbox{$\rhd$}\cr}}}
\def \coc {{^{\frown}}}
\def \lr {{\langle \rangle}}
\def \lor {{\langle 0 \rangle}}

\def \res {\lceil}

\def \la {\langle}
\def \ra {\rangle}
\def \Tau {\mathcal{T}}

\def \tp {\textrm{tp}}

\def \L {\mathcal{L}}

\def \M {\mathbb{M}}

\def\tp{\operatorname{tp}}

\newcommand{\uindep}[1][]{%
  \mathrel{
    \mathop{
      \vcenter{
        \hbox{\oalign{\noalign{\kern-3pt}\hfill\hspace{3.6pt}$\vert^u$\cr
              \noalign{\kern-.7ex}
              $\smile$\cr\noalign{\kern-.3ex}}}
      }
    }\displaylimits_{#1}
  }
}

\newcommand{\findep}[1][]{%
  \mathrel{
    \mathop{
      \vcenter{
        \hbox{\oalign{\noalign{\kern-3pt}\hfill\hspace{3.6pt}$\vert^f$\cr
              \noalign{\kern-.7ex}
              $\smile$\cr\noalign{\kern-.3ex}}}
      }
    }\displaylimits_{#1}
  }
}

\newcommand{\indep}[1][]{%
  \mathrel{
    \mathop{
      \vcenter{
        \hbox{\oalign{\noalign{\kern-.3ex}\hfil$\vert$\hfil\cr
              \noalign{\kern-.7ex}
              $\smile$\cr\noalign{\kern-.3ex}}}
      }
    }\displaylimits_{#1}
  }
}

\newcommand{\starindep}[1][]{%
  \mathrel{
    \mathop{
      \vcenter{
        \hbox{\oalign{\noalign{\kern-3pt}\hfill\hspace{3.5pt}$\vert^*$\cr
              \noalign{\kern-.7ex}
              $\smile$\cr\noalign{\kern-.3ex}}}
      }
    }\displaylimits_{#1}
  }
}

\newcommand{\ind}[1][]{%
  \mathrel{
    \mathop{
      \vcenter{
        \hbox{\oalign{\noalign{\kern-.3ex}\hfil$\vert$\hfil\cr
              \noalign{\kern-.7ex}
              $\smile$\cr\noalign{\kern-.3ex}}}
      }
    }\displaylimits_{#1}
  }
}

\newcommand{\inds}[1][]{%
  \mathrel{
    \mathop{
      \vcenter{
        \hbox{\oalign{\noalign{\kern-.3ex}\hfil$\shortmid$\hfil\cr
              \noalign{\kern-.7ex}
              $\smile$\cr\noalign{\kern-.3ex}}}
      }
    }\displaylimits_{#1}
  }
}

\def\notind#1#2{#1\setbox0=\hbox{$#1x$}\kern\wd0
\hbox to 0pt{\mathchardef\nn=12854\hss$#1\nn$\kern1.4\wd0\hss}
\hbox to 0pt{\hss$#1\mid$\hss}\lower.9\ht0 \hbox to 0pt{\hss$#1\smile$\hss}\kern\wd0}

\def\notstarind#1#2{#1\setbox0=\hbox{$#1x$}\kern\wd0
\hbox to 0pt{\mathchardef\nn=12854\hss$#1\nn$\kern1.4\wd0\hss}
\hbox to 0pt{\hss$#1\hspace{4.3pt}\mid^*$\hss}\lower.9\ht0 \hbox to 0pt{\hss$#1\smile$\hss}\kern\wd0}

\def\notfind#1#2{#1\setbox0=\hbox{$#1x$}\kern\wd0
\hbox to 0pt{\mathchardef\nn=12854\hss$#1\nn$\kern1.4\wd0\hss}
\hbox to 0pt{\hss$#1\hspace{5.2pt}\mid^f$\hss}\lower.9\ht0 \hbox to 0pt{\hss$#1\smile$\hss}\kern\wd0}

\title{Some Remarks on Kim-dividing in NATP Theories}

\author[J. Kim]{Joonhee Kim} 
\address{School of Mathematics\\
 Korea Institute for Advanced Study\\
 Seoul, South Korea}
 \email{kimjoonhee@kias.re.kr}

\author[H. Lee]{Hyoyoon Lee}
\address{Center for Nano Materials\\
 G-LAMP\\
 Sogang University\\
 Seoul, South Korea \newline
 \indent Department of Mathematics\\
 Yonsei University\\
 Seoul, South Korea}
 \email{hyoyoonlee@sogang.ac.kr, hyoyoonlee@yonsei.ac.kr}

\begin{document}

\maketitle

\begin{abstract}

In this note, we prove that Kim-dividing over models is always witnessed by a coheir Morley sequence, whenever the theory is NATP.

Following the strategy of Chernikov and Kaplan \cite{CK12}, we obtain some corollaries which hold in NATP theories. Namely, (i) if a formula Kim-forks over a model, then it quasi-divides over the same model, (ii) for any tuple of parameters $b$ and a model $M$, there exists a global coheir $p$ containing $\text{tp}(b/M)$ such that $B \ind^K_M b'$ for all $b'\models p|_{MB}$.

We also show that for coheirs in NATP theories, condition (ii) above is a necessary condition for being a witness of Kim-dividing, assuming that a witness of Kim-dividing exists (see Definition \ref{def:witness} in this note). That is, if we assume that a witness of Kim-dividing always exists over any given model, then a coheir $p\supseteq \text{tp}(a/M)$ must satisfy (ii) whenever it is a witness of Kim-dividing of $a$ over a model $M$. We also give a sufficient condition for the existence of a witness of Kim-dividing in terms of pre-independence relations.

At the end of the paper, we leave a short remark on Mutchnik's recent work \cite{Mut22}. We point out that the class of N-$\omega$-DCTP$_2$ theories, a subclass of the class of NATP theories, contains all NTP$_2$ theories and NSOP$_1$ theories. We also note that Kim-forking and Kim-dividing are equivalent over models in N-$\omega$-NDCTP$_2$ theories, where Kim-dividing is defined with respect to invariant Morley sequences, instead of coheir Morley sequences as in \cite{Mut22}.

\end{abstract}

\vspace{10pt}

\section{Introduction}

\subsection{Background}
Determining whether forking and dividing are equivalent seems to be a natural step in a study of a new dividing line. As in the study of simple theories, NTP$_2$ theories, and as in the recent successful study of NSOP$_1$ theories, the equivalence of forking and dividing serves as a starting point for many combinatorial model-theoretic observations. In this sense, the authors attempted to investigate the concept of Kim-forking and Kim-dividing in NATP theories and were able to obtain some partial results. This note is a report on them.

As a common extension of the class of NTP$_2$ theories and NSOP$_1$ theories (or NTP$_1$ theories, regarding the recent result in \cite{Mut22}), the class of NATP theories is expected to inherit some properties of these two classes. However, these features actually serve to reduce the number of things we can expect from NATP theories. That is, it tells us which properties we should not expect for forking and dividing in NATP theories. For example, in NATP theories, forking is not equivalent to dividing in general, even over a model, since there is an NSOP$_1$ theory where forking does not imply dividing \cite[Proposition 9.17]{KR20}.

 Thus it seems that other notions stronger than forking and dividing are needed. Fortunately, in \cite{KR20}, Kaplan and Ramsey introduce notions of Kim-forking and Kim-dividing, which are stronger than forking and dividing. `Stronger' here means that the notions are witnessed by `rarer' objects. Namely, dividing of a given formula is witnessed by an indiscernible sequence, while Kim-dividing of a given formula is witnessed by an invariant Morley sequence. We can extend this point of view. One can define a notion stronger than Kim-dividing using something rarer than invariant Morley sequences such as coheir Morley sequences, or strictly invariant Morley sequences introduced in \cite{CK12}.
 
 On the other hand, in the study of Kim-dividing in NATP theories, we cannot expect that every invariant Morley sequence over a given model will witness every instance of Kim-dividing over the model as in NSOP$_1$ theories. In fact, even coheir Morley sequences may not witness Kim-dividing of some formula in some NTP$_2$ theories. Let $M$ be a small model of DLO, $a,b$ live in the same cut of $M$ ({\it i.e.} $a\equiv_Mb$ and $m<a<m'$ for some $m,m'\in M$), and $a<b$.
Then the formula $\varphi(x,a,b):=a<x<b$ Kim-divides over $M$. If we choose any sequence $(a_ib_i)_{i<\omega}$ in the same cut with $a_0=a$, $b_0=b$ such that \[\cdots<a_{i+1}<a_i<\cdots<a_1<a_0<b_0<b_1<\cdots<b_i<b_{i+1}<\cdots,\]
then the sequence is a coheir Morley sequence over $M$ but $\{\varphi(x,a_i,b_i)\}_{i<\omega}$ is consistent.
 
 Thus in a study of NATP theories, if we want to build a class of indiscernible sequences that works meaningfully together with Kim's lemma, then the class should not include all coheir Morley sequences.

In order to show the equivalence between forking and dividing, the first thing we usually consider is establishing Kim's lemma. The first appearance of the statement is that:
\begin{fact}\label{fact:kim's lemma}\cite{Kim98}
In simple theories, if a formula $\varphi(x,a)$ divides over a set $A$, then $\{\varphi(x,a_i)\}_{i<\omega}$ is inconsistent for any Morley sequence $(a_i)_{i<\omega}$  over $A$ with $a_0=a$.
\end{fact}

Following this, Chernikov and Kaplan, Kaplan and Ramsey, and Mutchnik proved Kim's lemma for NTP$_2$, NSOP$_1$, and NSOP$_2$ theories, respectively, as follows.

\begin{itemize}
\item[\cite{CK12}] In NTP$_2$ theories, if a formula $\varphi(x,a)$ divides over a model $M$, then the partial type $\{\varphi(x,a_i)\}_{i<\omega}$ is inconsistent for any strictly invariant Morley sequence $(a_i)_{i<\omega}$ over $M$ with $a_0=a$.
\item[{\cite{KR20}}] In NSOP$_1$ theories, if a formula $\varphi(x,a)$ Kim-divides over a model $M$, then $\{\varphi(x,a_i)\}_{i<\omega}$ is inconsistent for any invariant Morley sequence $(a_i)_{i<\omega}$ over $M$ with $a_0=a$.
\item[{\cite{Mut22}}]  In NSOP$_2$ theories, if a formula $\varphi(x,a)$ canonical coheir-divides over a model $M$, 
then $\{\varphi(x,a_i)\}_{i<\omega}$ is inconsistent for any coheir Morley sequence $(a_i)_{i<\omega}$ over $M$ with $a_0=a$. (For the definition of canonical coheir-dividing, see Definition \ref{def:canonical coheir morley sequences})
\end{itemize}
We can see that the statements above are given in a uniform way, namely,

\begin{itemize}
\item[{[$\ast$]}] In $X$ theories, if a formula $Y$-divides, then the $Y$-dividing of the formula is witnessed by every $Z$ Morley sequence. 
\end{itemize}

The authors' aim is finding a suitable class of Morley sequences for $Z$ in the statement [$\ast$] above, assuming $X$=NATP and $Y$=Kim(=invariant). We could not finish our goal but found some partial results. As in the studies mentioned above and suggested in \cite{Kim09}, finding the condition $Z$ for Morley sequences is expected to play a crucial role in the study of independence in NATP.

\subsection{Overview}
First we observe that if we set $X$=NATP in the statement [$\ast$] above, then we can replace `Kim' with `coheir' for the variable $Y$.
Namely,
\begin{thm1}
Suppose that $T$ is NATP and let $M$ be a model of $T$. If $\varphi(x,a)$ Kim-divides over $M$, then it coheir-divides over $M$.
\end{thm1}

And following the strategy of Chernikov and Kaplan in \cite{CK12}, we obtain some corollaries.

\begin{cor2}
Suppose that the theory is NATP. If a formula Kim-forks over a model $M$, then it quasi-divides over $M$.
\end{cor2}

\begin{cor3}
Suppose that the theory is NATP. Then for each model $M$ and a tuple of parameters $b$, there exists a global coheir $p(x)$ over $M$ containing $\tp(b/M)$ such that $B\ind^K_M b'$ for all $b'\models p(x)|_{MB}$.
\end{cor3}

Although they still do not complete our goal, we also obtain both sufficient condition and necessary condition for $Z$ in [$\ast$], as follows.

\begin{prop4}
Let $T$ be an NATP theory and $\M$ its monster model. Suppose that for any model $M$ and $a\in\M$, there exists an invariant Morley sequence which is a witness of Kim-dividing of $a$ over $M$. Let $M$ be a model, $a\in\M$, and $I=(a_i)_{i<\omega}$ a coheir Morley sequence over $M$ with $a_0=a$. If $I$ is a witness of Kim-dividing of $a$ over $M$, then it is a strict coheir Morley sequence of $\tp(a/M)$.
\end{prop4}

So in particular, if there exists a class of coheir Morley sequences, say $Z$, which makes [$\ast$] hold with $X$=NATP and $Y$=Kim, then every $Z$-coheir Morley sequence is a strict coheir Morley sequence.

\begin{prop5}
Let $T$ be an NATP theory and $M$ a model. Then (i) implies (ii) where:
\begin{itemize}
\item[(i)] There exists a pre-independence relation $\ind$ which is stronger than $\ind^h$ and satisfies monotonicity, full existence, right chain condition for coheir Morley sequences, and strong right transitivity over $M$.
\item[(ii)] For all $b$, there exists a global coheir which is a witness of Kim-dividing of $b$ over $M$. Thus Kim-forking is equivalent with Kim-dividing over $M$.
\end{itemize}
\end{prop5}

More precisely, in the proof of Proposition \ref{prop:sufficient cond}, we show that if a pre-independence relation $\ind$ as above exists, then every global coheir $p$ generated by $\ind$ is always a witness of Kim-dividing of $b$ over $M$, where $b\models p|_M$. Note that if such a pre-independence relation $\ind$ exists, then the global coheir $p$ as above must exist by full existence of $\ind$.

\smallskip

 To find $Z$ in [$\ast$] under the assumption that $X$ is NATP, one of the ways we can try is constructing an antichain tree using two global types $p$ and $q$ corresponding to the properties indicated by $Y$=coheir and $Z$, respectively. The sequences generated by $p$ should be `path' parts of the antichain tree, and the sequences generated by $q$ should be `antichain' parts of the tree. We should make every antichain in the tree to be indices of a sequence generated by $q$, but there are infinitely many different types of antichains in tree structures with respect to $\L_0$. Thus we need a pre-independence relation strong enough to generate a global type with `good' properties that make its Morley sequences cover all antichain parts in a tree structure. For more detail on the role of the pre-independence relations on the construction of antichain trees, see the proofs of Theorem \ref{thm:Kim-dividing->coheir-dividing} and Proposition \ref{prop:sufficient cond}.

In his work \cite{Mut22}, Mutchnik introduces a new pre-independence relation $\ind^{\rm CK}$ and a tree property which is called $k$-DCTP$_2$ defined by using $2^{<\omega}$, whose consistent parts are descending combs, and inconsistent parts are paths. He proves that we can construct $k$-DCTP$_2$ by using a coheir $p$, and a canonical coheir $q$ whose existence is given by $\ind^{\rm CK}$. As a consequence, he shows that:

\begin{fact}\cite[Theorem 4.9]{Mut22}
If a theory does not have $k$-DCTP$_2$ for all $k<\omega$, and if a formula $\varphi(x,a)$ coheir-divides over a model $M$, then $\{\varphi(x,a_i)\}_{i<\omega}$ is inconsistent for any canonical coheir Morley sequence $(a_i)_{i<\omega}$ over $M$ with $a_0=a$.
\end{fact}

Thus in a theory which does not have $k$-DCTP$_2$ for all $k<\omega$ (we will call it an N-$\omega$-DCTP$_2$ theory), coheir-forking and coheir-dividing are equivalent over models. By combining our result Theorem \ref{thm:Kim-dividing->coheir-dividing} together, we can say that Kim-forking and Kim-dividing are equivalent over models in N-$\omega$-DCTP$_2$ theories. Interestingly, $k$-DCTP$_2$ implies SOP$_1$ and TP$_2$. And ATP implies $k$-DCTP$_2$ for any $k<\omega$. So we get a new dividing line inside NATP, outside both NTP$_1$ and NTP$_2$, where Kim-forking and Kim-dividing are equivalent over models. We will discuss this in Section 5.

\section{Preliminary}
In this section we enumerate the basic notions and facts on the subject we will cover throughout this paper. More basic notions on model theory that we do not mention in this section follow \cite{TZ}.

\subsection{Tree properties}
Let us recall basic notations and facts about tree properties.

\begin{notation}\label{language of trees}
Let $\kappa$ and $\lambda$ be cardinals.
\begin{enumerate}
\item[(i)] By $\lambda^\kappa$, we mean the set of all functions from $\lambda$ to $\kappa$. 
\item[(ii)] By $\trec$, we mean $\bigcup_{\alpha<\kappa}{\lambda^\alpha}$ and call it a tree. If $\lambda=2$, we call it a binary tree. If $\lambda\geq\omega$, then we call it an infinitary tree.
\item[(iii)] By $\emptyset$ or $\lr$, we mean the empty string in $\trec$, which means the empty set (recall that every function can be regarded as a set of ordered pairs).
\end{enumerate}
Let $\eta,\nu\in {\trec}$. 
\begin{enumerate}
\item[(iv)] By $\eta\tri\nu$, we mean $\eta \subseteq \nu$. So $\trec$ is partially ordered by $\tri$. If $\eta\tri\nu$ or $\nu\tri\eta$, then we say $\eta$ and $\nu$ are comparable. 
\item[(v)] By $\eta\perp\nu$, we mean that $\eta\not\tri\nu$ and $\nu\not\tri\eta$. We say $\eta$ and $\nu$ are incomparable if $\eta\perp\nu$.
\item[(vi)] By $\eta\wedge\nu$, we mean the maximal $\xi\in{\trec}$ such that $\xi\tri\eta$ and $\xi\tri\nu$.
\item[(vii)] By $l(\eta)$, we mean the domain of $\eta$.
\item[(viii)] By $\eta\lex\nu$, we mean that either $\eta\tri\nu$, or $\eta\perp\nu$ and $\eta(l(\eta\wedge\nu))<\nu(l(\eta\wedge\nu))$. 
\item[(ix)] By $\eta\coc\nu$, we mean $\eta\cup\{(l(\eta)+i,\nu(i)):i< l(\nu)\}$. Note that $\emptyset\coc\nu$ is just $\nu$.
\end{enumerate}
Let $X\subseteq \trec$.
\begin{enumerate}
\item[(x)]
By $\eta\coc X$ and $X\coc\eta$, we mean $\{\eta\coc x:x\in X\}$ and $\{x\coc\eta:x\in X\}$ respectively.
\end{enumerate}
Let $\eta_0,...,\eta_n\in\trec$.
\begin{enumerate}
\item[(xi)] We say a subset $X$ of $\trec$ is an {\it antichain} if the elements of $X$ are pairwise incomparable, {\it i.e.}, $\eta\perp\nu$ for all $\eta,\nu\in X$ with $\eta\neq\nu$.
\end{enumerate}
\end{notation}

\begin{dfn}\cite{DS04}\cite{AK20}\cite{Che14}\cite{CR16}
Let $\varphi(x,y)$ be an $\mathcal{L}$-formula. 
\begin{enumerate}
\item[(i)] We say $\varphi(x,y)$ has the {\it tree property} (TP) if there exists a tree-indexed set $(a_\eta)_{\eta\in\treo}$ of parameters and $k\in\omega$ such that
\begin{enumerate}
\item[] $\{\varphi(x,a_{\eta_{\res n}})\}_{n\in\omega}$ is consistent for all $\eta\in\troo$ (path consistency), and
\item[] $\{\varphi(x,a_{\eta\coc i})\}_{i\in\omega}$ is $k$-inconsistent for all $\eta\in\treo$, {\it i.e.}, any subset of $\{\varphi(x,a_{\eta\coc i})\}_{i\in\omega}$ of size $k$ is inconsistent.
\end{enumerate}
\item[(ii)] We say $\varphi(x,y)$ has the {\it tree property of the first kind} (TP$_1$) if there is a tree-indexed set $(a_\eta)_{\eta\in\treo}$ of parameters such that
\begin{enumerate}
\item[] $\{\varphi(x,a_{\eta_{\res n}})\}_{n\in\omega}$ is consistent for all $\eta\in\troo$, and
\item[] $\{\varphi(x,a_\eta),\varphi(x,a_\nu)\}$ is inconsistent for all $\eta\perp\nu\in\treo$.
\end{enumerate}
\item[(iii)] We say $\varphi(x,y)$ has the {\it tree property of the second kind} (TP$_2$) if there is an array-indexed set $(a_{i,j})_{i,j\in\omega}$ of parameters such that
\begin{enumerate}
\item[] $\{\varphi(x,a_{n,\eta(n)})\}_{n\in\omega}$ is consistent for all $\eta\in\troo$, and
\item[] $\{\varphi(x,a_{i,j}),\varphi(x,a_{i,k}))\}$ is inconsistent for all $i,j,k\in\omega$ with $j\neq k$.
\end{enumerate}
\item[(iv)] We say $\varphi(x,y)$ has the {\it 1-strong order property} (SOP$_1$) if there is a binary-tree-indexed set $(a_\eta)_{\eta\in\tree}$ of parameters such that
\begin{enumerate}
\item[] $\{\varphi(x,a_{\eta_{\res n}})\}_{n\in\omega}$ is consistent for all $\eta\in2^{\omega}$,
\item[] $\{\varphi(x,a_{\eta\coc 1}),\varphi(x,a_{\eta\coc 0\coc\nu})\}$ is inconsistent for all $\eta,\nu\in\tree$. 
\end{enumerate}
\item[(v)] We say $\varphi(x,y)$ has the {\it 2-strong order property} (SOP$_2$) if there is a binary-tree-indexed set $(a_\eta)_{\eta\in\tree}$ of parameters such that
\begin{enumerate}
\item[] $\{\varphi(x,a_{\eta_{\res n}})\}_{n\in\omega}$ is consistent for all $\eta\in2^\omega$,
\item[] $\{\varphi(x,a_\eta),\varphi(x,a_\nu)\}$ is inconsistent for all $\eta\perp\nu\in\tree$.
\end{enumerate}
\item[(vi)] We say $\varphi(x,y)$ has the {\it antichain tree property} (ATP) if there is a binary-tree-indexed set $(a_\eta)_{\eta\in\tree}$ of parameters such that
\begin{enumerate}
\item[] $\{\varphi(x,a_\eta)\}_{\eta\in X}$ is consistent for all antichain $X \subseteq2^{<\omega}$,
\item[] $\{\varphi(x,a_\eta),\varphi(x,a_\nu)\}$ is inconsistent for all $\eta,\nu\in\tree$ with $\eta\trn\nu$.
\end{enumerate}
\item[(vii)] We say a theory has TP (or is TP) and call it a TP theory if there is a formula having TP with respect to its monster model. We define TP$_1$, TP$_2$, SOP$_1$, SOP$_2$, and ATP theories in the same manner.
\item[(viii)] We say a theory is NTP if the theory is not TP. We define NTP$_1$, NTP$_2$, NSOP$_1$, NSOP$_2$, and NATP theories in the same manner.
\end{enumerate}
\end{dfn}

\begin{remark}
By compactness, we can replace $2^{<\omega}$ in the definitions of SOP$_2$ and ATP with $\trec$ for any cardinal $\lambda$, and an infinite cardinal $\kappa$.
\end{remark}

The following facts are from \cite{Con12}, \cite{DS04}, \cite{KK11}, \cite{She90}, \cite{AK20}, and \cite{Mut22}.

\begin{fact}
\begin{enumerate}
\item[(i)] A theory has TP$_1$ if and only if it has SOP$_2$. 
\item[(ii)] A theory has TP if and only if it has TP$_1$ or TP$_2$.
\item[(iii)] A theory has SOP$_2$ if and only if it has SOP$_1$.
\item[(iv)] If a theory has ATP, then it has SOP$_1$ and TP$_2$.
\end{enumerate}
\end{fact}

Thus we have the following diagram, 

\begin{center}
\includegraphics[width=0.53\linewidth]{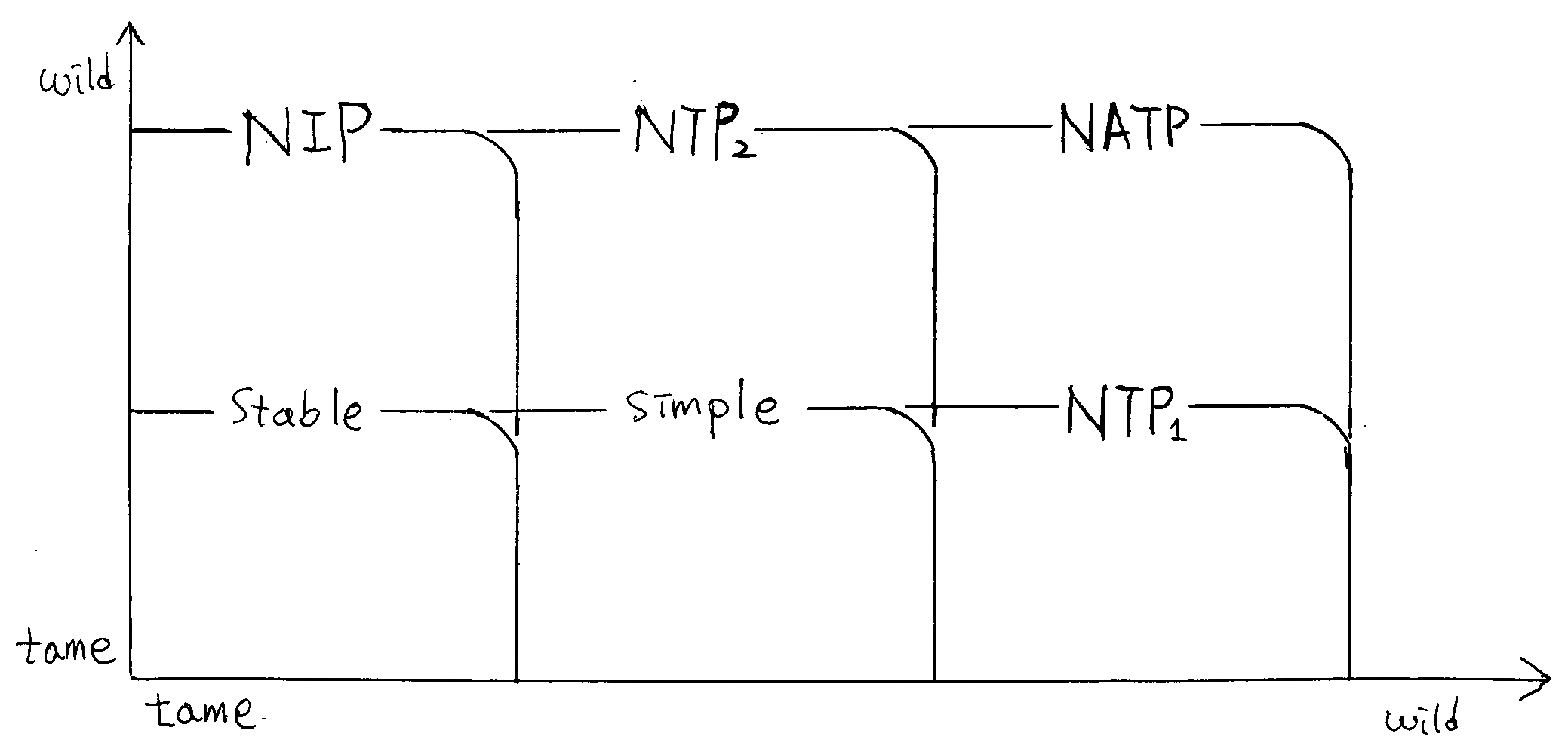}
\end{center}
where simplicity is equivalent to not having TP.



\subsection{Tree indiscernibility and modeling property}

Now we recall the notion and some facts on tree indiscernibility and modeling property. Proof for the most facts can be found in \cite{KK11}, \cite{TT12}, \cite{KKS14}, and \cite{Sco15}. Let $\mathcal{L}_0=\{\tri,\lex,\wedge\}$ be a language where $\tri$, $\lex$ are binary relation symbols, and $\wedge$ is a binary function symbol. Then for cardinals $\kappa>1$ and $\lambda$, a tree $\trec$ can be regarded as an $\mathcal{L}_0$-structure whose interpretations of $\tri,\lex,\wedge$ follow Notation \ref{language of trees}.

\begin{dfn}
Let $\overline{\eta}=( \eta_0,...,\eta_n)$ and $\overline{\nu}=( \nu_0,...,\nu_n)$ be finite tuples of $\trec$. 
\begin{enumerate}
\item[(i)] By $\textrm{qftp}_0(\overline{\eta})$, we mean the set of quantifier free $\mathcal{L}_0$-formulas $\varphi(\overline{x})$ such that $\trec\models \varphi(\overline{\eta})$. 
\item[(ii)] By $\overline{\eta}\sim_0\overline{\nu}$, we mean $\textrm{qftp}_0(\overline{\eta})=\textrm{qftp}_0(\overline{\nu})$. We say $\bar{\eta}$ and $\bar{\nu}$ are {\it strongly isomorphic} if $\overline{\eta}\sim_0\overline{\nu}$.
\end{enumerate}
Let $\mathcal{L}$ be a language, $T$ be a complete $\mathcal{L}$-theory, $\mathbb{M}$ be a monster model of $T$, and $( a_\eta )_{\eta\in\trec}$, $( b_\eta )_{\eta\in\trec}$ be tree-indexed sets of parameters from $\mathbb{M}$.
\begin{enumerate}
\item[(iii)] We say $( a_\eta )_{\eta\in\trec}$ is {\it strongly indiscernible over $A$} if $\overline{a}_{\overline{\eta}}\equiv_A\overline{b}_{\overline{\nu}}$ for all $\textrm{qftp}_0(\overline{\eta})=\textrm{qftp}_0(\overline{\nu})$. If $A=\emptyset$, then we just say it is strongly indiscernible.
\item[(iv)] We say $( b_\eta)_{\eta\in\trec}$ is {\it strongly locally based on $(a_\eta)_{\eta\in\trec}$ over $A$} if for all $\overline{\eta}=(\eta_0,...,\eta_n)$ and a finite set of $\mathcal{L}(A)$-formulas $\Delta$, there is $\overline{\nu}=(\nu_0,...,\nu_n)$ such that $\overline{\eta}\sim_0\overline{\nu}$ and $\textrm{tp}_\Delta(\overline{a}_{\overline{\eta}})=\textrm{tp}_\Delta(\overline{b}_{\overline{\nu}})$, where $\overline{a}_{\overline{\eta}}$ and $\overline{b}_{\overline{\nu}}$ denote $( a_{\eta_0},...,a_{\eta_n})$ and $(b_{\nu_0},...,b_{\nu_n})$ respectively. If $A=\emptyset$, then we just say it is strongly locally based on $(a_\eta)_{\eta\in\trec}$.
\end{enumerate}
\end{dfn}

\begin{fact}\label{modeling property}
Let $\L$ be a language and $\M$ a sufficiently saturated $\L$-structure, and $( a_\eta )_{\eta\in\treo}$ a tree-indexed set of parameters in $\M$. Then there is a strongly indiscernible $( b_\eta )_{\eta\in\treo}\subseteq\M$ which is strongly locally based on $( a_\eta )_{\eta\in\treo}$. 
\end{fact}

\noindent The proof can be found in \cite{KK11}, \cite{TT12}, and \cite{Sco15}. 
Note that in some context, for a given set $A$, we may assume $( b_\eta )_{\eta\in\treo}$ is strongly indiscernible over $A$ and strongly locally based on $( a_\eta )_{\eta\in\treo}$ over $A$, by adding constant symbols to $\L$.

The above statement is called {\it the modeling property of strong indiscernibility} (in short, we write it {\it the strong modeling property}). More precisely, we say an $\L_0$-structure $I$ has {\it the strong modeling property} if there always exists strongly indiscernible $(b_i)_{i\in I}$ which is strongly locally based on $(a_i)_{i\in I}$, for any given $(a_i)_{i\in I}$.

\subsection{Pre-independence relations, invariant types, and forking}

We quote the following notions of pre-independence relations from \cite{Adl07}, \cite{Adl08}, and \cite{CK12}.

\begin{dfn}\label{def:property_of_ind}\cite{Adl07}\cite{Adl08}\cite{Adl09}\cite[Definition 2.4]{CK12}
A {\it pre-independence relation} is an invariant ternary relation $\ind$ on sets.
 If a triple of sets $(a,b,c)$ is in the pre-independence relation $\ind$, then we write it $a\ind_c b$ and say {\it ``$a$ is $\ind$-independent from $b$ over $c$''}. Throughout this paper we will consider the following properties for a pre-independence relation. (If it is clear in the context, then we omit the words in the parentheses.)
\begin{itemize}
\item[(i)] Monotonicity (over $d$): If $aa'\ind_d bb'$, then $a\ind_d  b$.
\item[(ii)] Base monotonicity (over $d$): If $a\ind_d bb'$, then $a\ind_{db}b'$.
\item[(iii)] Left transitivity (over $d$): If $a\ind_{db} c$ and $b\ind_dc$, then $ab\ind_d c$.
\item[(iv)] Strong left transitivity (over $d$): If $a\ind_{d} bc$ and $b\ind_dc$, then $ab\ind_d c$.
\item[(v)] Strong right transitivity (over $d$): If $ab\ind_{d} c$ and $a\ind_db$, then $a\ind_d bc$.
\item[(vi)] Left extension (over $d$): If $a\ind_d b$, then for all $c$, there exists $c'\equiv_{da} c$ such that $ac'\ind_d b$. 
\item[(vii)] Right extension (over $d$): If $a\ind_d b$, then for all $c$, there exists $c'\equiv_{db} c$ such that $a\ind_d bc'$. 
\item[(viii)] Left existence (over $d$):
$a\ind_dd$ for all $a$. We say a set $d$ is an {\it extension base for $\ind$} if $\ind$ satisfies left existence over $d$.
\item[(ix)] Full existence (over $d$): For all $a,b$, there exists $a'\equiv_d a$ such that $a'\ind_db$. Equivalently, there exists $b'\equiv_d b$ such that $a\ind_d b'$.
\item[(x)] Finite character (over $d$): If $a\not\ind_db$, then there exist finite $a'\subseteq a$ and $b'\subseteq b$ such that $a'\not\ind_db'$. 
\item[(xi)] Strong finite character (over $d$): If $a\not\ind_db$, then there exist finite subtuple $b'\subseteq b$, finite tuples $x'$, $y'$ of variables with $|x'|\le|a|$, $|y'|=|b'|$, and a formula $\varphi(x',y')\in\L(d)$ such that $\varphi(x',b')\in\tp(a/db)$ and $a'\not\ind_db'$ for all $a'\models\varphi(x',b')$.
\end{itemize}
\end{dfn}

To consider more properties for pre-independence relations, we need the following definition.

\begin{dfn}
Let $M$ be a model. 
\begin{itemize}
\item[(i)] We say a complete type $p(x)$ is {\it invariant over $M$} if $\varphi(x,b)\leftrightarrow\varphi(x,b')\in p(x)$ for all $\varphi(x,y)$ and $b\equiv_M b '$.
\item[(ii)] We say a type $p(x)$ is {\it finitely satisfiable in $M$} or a {\it coheir over $M$} if for all finite subsets $\Delta(x)\subseteq p(x)$, there exists $m\in M$ such that $\models\Delta(m)$. 
\item[(iii)] We say a type $p(x)$ is an {\it heir over $M$} if for all $\varphi(x,b)\in p(x)$, there exists $m\in M$ such that $\varphi(x,m)\in p(x)$.
\end{itemize}
\end{dfn}

\begin{rmk}\label{rmk:the number of inv tp}
It is clear from the definition that if a global type is a coheir over $M$, then it is invariant over $M$ for any model $M$. Note that for any given model $M$ and a tuple of parameters $a$, and for any $|M|^+$-saturated model $N$ containing $M$, two global invariant types $p$ and $q$ over $M$ containing $\tp(a/M)$ are the same if and only if $p|_N=q|_N$ \cite[Section 2]{Sim}. Thus there exist only boundedly many invariant global types over $M$ containing $\tp(a/M)$. 
\end{rmk}

\begin{dfn}
Let $a$ be a tuple of parameters, $M$ a model, and $\kappa$ a cardinal.
\begin{itemize}
\item[(i)] A sequence $(a_i)_{i<\kappa}$ is called an {\it invariant Morley sequence in $\tp(a/M)$} if there exists a global invariant type $p(x)\supseteq\tp(a/M)$ over $M$ such that $a_i\models p(x)|_{Aa_{<i}}$ for all $i<\kappa$.
\item[(ii)] A sequence $(a_i)_{i<\kappa}$ is called a {\it coheir Morley sequence in $\tp(a/M)$} if there exists a global coheir $p(x)\supseteq\tp(a/M)$ over $M$ such that $a_i\models p(x)|_{Ma_{<i}}$ for all $i<\kappa$.
\end{itemize}
\end{dfn}

\begin{dfn}
Continuing Definition \ref{def:property_of_ind}, let $\ind$ be a pre-independence relation. We say it satisfies {\it the right chain condition (for coheir Morley sequences, over $d$)} if it satisfies the following condition. 
\begin{itemize}
\item[(xii)] If $a\ind_db$ and $(b_i)_{i<\omega}$ is a coheir Morley sequence over $M$ with $b_0=b$, then there exists $a'$ such that $a'\equiv_{db}a$, $a'\ind_d(b_i)_{i<\omega}$, and 
$(b_i)_{i<\omega}$ is indiscernible over $Ma'$.
\end{itemize}
We omit the words in the parentheses if it is clear in the context.
\end{dfn}

Now we recall the notions of forking and dividing.

\begin{dfn}\cite{Kim98}\cite{CK12}\cite{KR20}\label{def:dividing forking}
Let $\varphi(x,y)$ be a formula, $a$ a tuple of parameters with $|a|=|y|$, and $A$ a set.
\begin{itemize}
\item[(i)] We say {\it $\varphi(x,a)$ divides over $A$} if there exists an indiscernible sequence $(a_i)_{i<\omega}$ with $a_0=a$ such that $\{\varphi(x,a_i)\}_{i<\omega}$ is inconsistent.
\item[(ii)] We say {\it $\varphi(x,a)$ Kim-divides over $A$ if there exists an invariant Morley sequence $(a_i)_{i<\omega}$ in $\tp(a/A)$ such that $\{\varphi(x,a_i)\}_{i<\omega}$ is inconsistent.}
\item[(iii)] We say {\it $\varphi(x,a)$ coheir-divides over $A$} if there exists a coheir Morley sequence $(a_i)_{i<\omega}$ in $\tp(a/A)$  such that $\{\varphi(x,a_i)\}_{i<\omega}$ is inconsistent.
\item[(iv)]We say {\it $\varphi(x,a)$ quasi-divides over $A$} if there exist $a_0,...,a_n$ such that $a_i\equiv_A a$ and $\{\varphi(x,a_i)\}_{i\le n}$ is inconsistent.
\item[(v)] We say a type $p(x)$ {\it forks over $A$} if there exist $\psi_0(x),...,\psi_n(x)$ with parameters such that $p(x)\vdash \bigvee_{i\le n} \psi_i(x)$ and $\psi_i(x)$ divides over $A$ for each $i\le n$.
\item[(vi)] We say a type $p(x)$ {\it Kim-forks over $A$} if there exist $\psi_0(x),...,\psi_n(x)$ with parameters such that $p(x)\vdash \bigvee_{i\le n} \psi_i(x)$ and $\psi_i(x)$ Kim-divides over $A$ for each $i\le n$.
\item[(vii)] A sequence $(a_i)_{i<\kappa}$ is called a {\it Morley sequence in $\tp(a/A)$} if it is indiscernible over $A$, $a_0=a$, and $\tp(a_i/Aa_{<i})$ does not fork over $A$.
\end{itemize}
\end{dfn}

\begin{dfn}\label{def_independence}
Let $M$ be a model.
\begin{itemize}
\item[(i)] If there exists a global invariant type $p(x)\supseteq\tp(a/Mb)$ over $M$ then we write $a\ind^i_M b$.
\item[(ii)] If $\tp(a/Mb)$ is a coheir over $M$, then we write $a\ind^u_M b$.
\item[(iii)] If $\tp(a/Mb)$ is an heir over $M$, then we write $a\ind^h_M b$ .
\item[(iv)] If $\tp(a/Mb)$ does not fork over $M$, then we write $a\ind^f_M b$.
\item[(v)] If $\tp(a/Mb)$ does not Kim-fork over $M$, then we write $a\ind^K_M b$.
\item[(vi)] We write $a\ind^{d}_M b$ if $\tp(a/Mb)$ has no formula dividing over $M$.
\item[(vii)] We write $a\ind^{K\!d}_M b$ if $\tp(a/Mb)$ has no formula Kim-dividing over $M$.
\item[(viii)] We write $a\ind^{cd}_M b$ if $\tp(a/Mb)$ has no formula coheir-dividing over $M$.
\end{itemize}
\end{dfn}

\begin{fact}\label{fact:i,u,f,K-ind}\cite{Adl08}\cite{CK12}
Let $M$ be a model and $a,b$ tuples of parameters.
\begin{itemize}
\item[(i)] $\ind^u$, $\ind^i$, and $\ind^f$ are pre-independence relations and satisfy monotonicity, base monotonicity, finite character, strong finite character, left transitivity, and right extension over $M$.
\item[(ii)] Additionally, $\ind^u$ satisfies left extension over $M$.
\item[(iii)] By base monotonicity and left transitivity, $\ind^i$, $\ind^u$, and $\ind^f$ satisfy strong left transitivity over $M$.
\item[(v)] $a\ind^u_M b$ implies $a\ind^i_M b$, $a\ind^i_M b$ implies $a\ind_M^f b$, and $a\ind_M^f b$ implies $a\ind_M^K b$.
\end{itemize}
\end{fact}

\section{Kim-dividing and coheir-dividing in NATP theories}

 The main statement of this section is that: if a formula $\varphi(x,a)$ Kim-divides over a model $M$, then it coheir-divides over $M$. From this result, we will observe some corollaries in a similar way to what Chernikov and Kaplan do in their work on NTP$_2$ theories \cite{CK12}.

 First we recall the notion of ill-founded trees which is introduced by Kaplan and Ramsey \cite{KR20}. The following definition is slightly different from the original \cite[Definition 5.2]{KR20}, but the idea is the same.

\begin{dfn}\label{def_Tau}
Suppose $\alpha$ and $\delta$ are ordinals. We define $\Tau_\alpha^\delta$ to be the set of functions $\eta$ so that
\begin{enumerate}
\item[(i)] dom($\eta$) is an end-segment of $\alpha$ of the form $[\beta,\alpha)$ for $\beta$ equal to $0$ or a successor ordinal. If $\alpha$ is a successor or $0$, we allow $\beta=\alpha$, {\it i.e.} dom($\eta$)=$\emptyset$. Note that $\Tau_0^\delta=\{\emptyset\}$.
\item[(ii)] ran($\eta$)$\subseteq\delta$.
\item[(iii)] Finite support: the set $\{\gamma\in{\rm dom}(\eta):\eta(\gamma)\neq0\}$ is finite.
\end{enumerate}
 Let $\L_1=\{\tri,\lex,\wedge,\len\}$ and $\L_{s,\alpha}=\{\tri,\lex,\wedge,\{P_\beta\}_{\beta<\alpha}\}$ for each ordinal $\alpha$, where $\len$ is a binary relation symbol, $P_\beta$ is an unary relation symbol for each $\beta$. 
We interpret $\Tau_\alpha^\delta$ as an $\L_0$-structure, $\L_1$-structure and an $\L_{s,\alpha}$-structure by defining each symbol as below.
\begin{enumerate}
\item[(iv)] $\eta\tri \nu$ if and only if $\eta\subseteq\nu$. Write $\eta\perp\nu$ if $\neg(\eta\tri\nu)$ and $\neg(\nu\tri \eta)$.
\item[(v)] $\eta\wedge\nu=\eta|_{[\beta,\alpha)}=\nu|_{[\beta,\alpha)}$ where $\beta={\rm min}\{\gamma : \eta|_{[\gamma,\alpha)}=\nu|_{[\gamma,\alpha)}\}$, if non-empty (note that $\beta$ will not be a limit, by finite support). Define $\eta\wedge\nu$ to be the empty function if this set is empty (note that this cannot occur if $\alpha$ is a limit).
\item[(vi)] $\eta\lex\nu$ if and only if $\eta\lhd\nu$ or, $\eta\perp\nu$ with dom($\eta\wedge\nu$)=$[\gamma+1,\alpha)$ and $\eta(\gamma)<\nu(\gamma)$.
\item[(vii)] For each ordinal $\beta<\alpha$, let $P_\beta^{\alpha,\delta}=\{\eta\in\Tau_\alpha^\delta:{\rm dom}(\eta)=[\beta,\alpha)\}$ (the $\beta$-th level in $\Tau_\alpha^\delta$). 
If it is clear in the context, we omit $\alpha$ and $\delta$, just write $P_\beta$. Note that if $\beta$ is limit then $P_\beta$ is empty.
\item[(viii)] $\eta\len\nu$ if ${\rm dom}(\eta)\supsetneq{\rm dom}(\nu)$.
\end{enumerate}
We will also need the following notation.
\begin{enumerate}
\item[(ix)] Canonical inclusion: For $\alpha<\alpha'$, $\Tau_\alpha^2$ can be embedded in $\Tau_{\alpha'}^2$ with respect to $\L_{s,\alpha'}$ by a map $f_{\alpha,\alpha'}:\Tau_\alpha^2\to \Tau_{\alpha'}^2:\eta\mapsto \eta\cup\{(\beta,0):\beta\in\alpha'\setminus\alpha\}$. Unless otherwise stated, we regard $\Tau_\alpha^2$ as $f_{\alpha,\alpha'}(\Tau_\alpha^2)$ in $\Tau_{\alpha'}^2$. Note that by finite support, $\Tau_\alpha^2$ can be regarded as $\bigcup_{\beta<\alpha}\Tau_\beta^2$ with respect to canonical inclusions, for each limit ordinal $\alpha$.
\item[(x)] $\eta\perp_{lex}\nu$ if and only if $\eta\lex\nu$ and $\eta\not\tri\nu$. For an indexed set $\{a_\eta\}_{\eta\in\Tau_\alpha^2}$ and $\eta\in\Tau_\alpha^2$, by $a_{\perp_{lex}\eta}$ we mean the set $\{a_\nu:\nu\perp_{lex}\eta\}$.
\item[(xi)] For each $\eta\in\Tau_\alpha^\delta$, let $t(\eta)$ be an ordinal such that ${\rm dom}(\eta)=[t(\eta),\alpha)$.
\item[(xii)] For each $\eta\in\Tau_\alpha^\delta$, let $C_\eta^{\alpha,\delta}$ (the cone on $\eta$ in $\Tau_\alpha^\delta$) be the set of all $\nu\in\Tau_\alpha^\delta$ such that $\eta\tri\nu$. 
If it is clear in the context, we omit $\alpha$ and $\delta$, just write $C_\eta$.
\item[(xiii)] For each $\eta\in\Tau_\alpha^\delta$, $i<\delta$, let $(i)\coc\eta=\eta\cup\{(\alpha,i)\}\in\Tau_{\alpha+1}^\delta$.
\item[(xiv)] $(i)\coc\Tau_\alpha^\delta=\{(i)\coc\eta:\eta\in\Tau_\alpha^\delta\}\subseteq\Tau_{\alpha+1}^\delta$ for each $\alpha$, $\delta$ and $i<\delta$.
\end{enumerate}
\end{dfn}

To simplify the argument of proof of the main statement, we need some definitions and notations.

\begin{rmk}\label{rmk:Tau}
Let $\alpha,\beta$ be ordinals and suppose $\beta<\alpha$. Then $\lex$ is a well-ordering on $P_\beta^{\alpha,2}$.
\begin{proof}
It is enough to show that $P_0^{\alpha,2}$ of $\Tau_\alpha^2$ is well-ordered by $\lex$. We show this using induction on $\alpha$. Clearly $P_0^{0,2}$ of $\Tau_0^2$ is well ordered by $\lex$. Suppose that for all $\beta<\alpha$, $P^{\beta,2}_0$ of $\Tau_{\beta}^2$ is well ordered by $\lex$. To get a contradiction, we assume $P^{\alpha,2}_0$ of $\Tau_\alpha^2$ is not well ordered by $\lex$. Then there exists a decreasing $(\eta_i)_{i<\omega}$ in $P^{\alpha,2}_0$ of $\Tau_\alpha^2$. By finite support, for each $i<\omega$, there exists $\beta_i<\alpha$ such that $\eta_i(\beta_i)=1$ and $\eta_i(\beta)=0$ for all $\beta_i<\beta<\alpha$. Note that $\beta_i\ge \beta_{i+1}$ for all $i<\omega$. Since $\alpha$ is well ordered by $<$, there exists $i^*<\omega$ such that $\beta_{i^*}=\beta_i$ for all  $i>i^*$. Let $\eta'_i:=\eta_i|_{\beta_{i^*}}$ for each $i>i^*$. Note that the domain of each $\eta_i$ is $\alpha$ since $\eta_i\in P^{\alpha,2}_0$. Thus $\eta'_i\in P^{\beta_{i^{*}}\!,2}_0$ of $\Tau_{\beta_{i^*}}^2$ for all $i>i^*$, and hence $(\eta'_i)_{i^*<i<\omega}$ forms a descending sequence in $P^{\beta_{i^{*}}\!,2}_0$ of $\Tau_{\beta_{i^*}}^2$, which yields a contradiction with the induction hypothesis.
\end{proof}
\end{rmk}

\begin{nota}
For ordinals $\beta<\alpha$, let $\Omega_\beta^{\alpha,2}$ be the ordinal which is order isomorphic to $(P_\beta^{\alpha,2},\lex)$ in $\Tau_\alpha^2$, and let 
$(\eta(\alpha,\beta,i))_{i<\Omega^{\alpha,2}_\beta}$ be the enumeration of $P_\beta^{\alpha,2}$ with respect to $\lex$.
\end{nota}
\begin{dfn}\cite[Section 2.2.1]{Sim}
Let $A$ be a set, $p$ and $q$ be global $A$-invariant types. By $p(x)\otimes q(y)$, we mean the global type with variables $x$ and $y$ such that for all $A\subseteq B\subseteq \M$ and $\L(B)$-formula $\varphi(x,y)$, $\varphi(x,y)\in p(x)\otimes q(y)$ if and only if there exists $b\models q|_{B}$ such that $\varphi(x,b)\in p(x)$.

Let $I$ be a linearly ordered set and $(q_i(x_i))_{i\in I}$ a sequence of global $A$-invariant types. By $\bigotimes_{i\in I}q_i(x_i)$, we mean the global type with variables $(x_i)_{i\in I}$ such that for all $A\subseteq B\subseteq \M$ and $\L(B)$-formula $\varphi(x_{i_0},...,x_{i_n})$ with $i_0<...<i_n$ and $i_0,...,i_n\in I$, $\varphi(x_{i_0},...,x_{i_n})\in \bigotimes_{i\in I}q_i(x_i)$ if and only if there exist $a_{i_0},...,a_{i_{n}}$ such that $\models \varphi(a_{i_0},...,a_{i_n})$ and $a_{i_k}\models q_{i_k}|_{Ba_{i_{<k}}}$ for all $k\le n$.
\end{dfn}

\begin{fact}\cite[Fact 2.19, 2.20]{Sim}\label{fact:type product}
If $(q_i)_{i\in I}$ is a sequence of $A$-invariant types, then $\bigotimes_{i\in I}q_i(x_i)$ is also an $A$-invariant type. Moreover, if the invariant types are coheirs, then the type product $\bigotimes_{i\in I}q_i(x_i)$ is also a coheir. Let $(I,<_I)$ and $(J,<_J)$ be linearly ordered sets. Let $I\oplus J$ be a disjoint union of $I$, $J$ and give a linear order $<$ on $I\oplus J$ by $x<y$ if and only if $x,y\in I$ and $x<_I y$, $x,y\in J$ and $x<_J y$, or $x\in I$ and $y\in J$. Then for given two sequences of global $A$-invariant types $(q_i)_{i\in I}$ and $(q_j)_{j\in J}$, $\left(\bigotimes_{j\in J}q_j(x_j)\right)\otimes \left(\bigotimes_{i\in I}q_i(x_i)\right)=\bigotimes_{k\in I\oplus J}q_k(x_k)$. In other words, the type product $\otimes$ on invariant types is associative.
\end{fact}

\begin{nota}\label{nota:variable rest}
Let $\alpha$ be an ordinal and $p((x_\eta)_{\eta\in\Tau_\alpha^2})$ a global type. For a subset $X$ of $\Tau_\alpha^2$, let $p|_X:=p|_{(x_\eta)_{\eta\in X}}$.
\end{nota}

\begin{dfn}
For each $\eta\in\Tau_\alpha^2$, $\{\xi_0,...,\xi_n\}$ with $\xi_0\lex\cdots\lex\xi_n$ in $\Tau_\alpha^2$ is called {\it the bottom antichain of $\eta$} if they form an antichain and satisfy $\{x_\nu\}_{\nu\perp_{lex}\eta}=\bigcup_{i\le n}C_{\xi_i}$. 
\end{dfn}

\begin{center}
\includegraphics[width=0.6\linewidth]{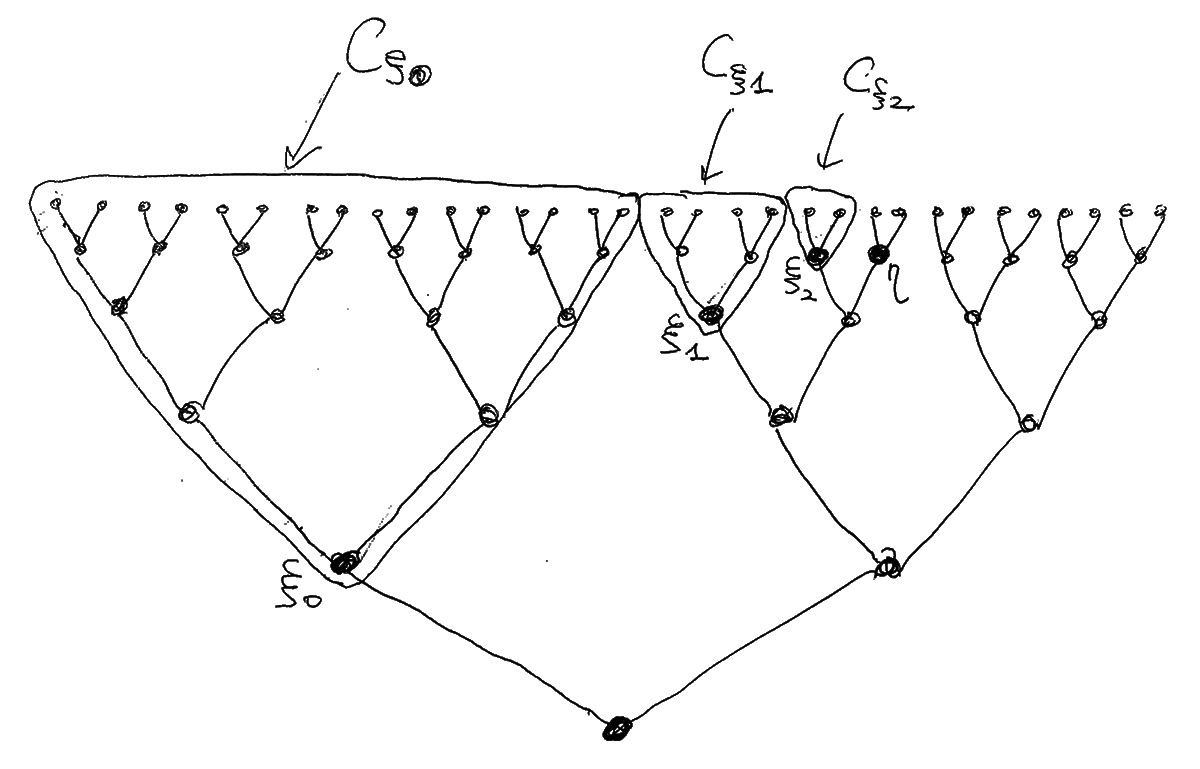}\\
\center {\scriptsize fig. The bottom antichain $\{\xi_0,\xi_1,\xi_2\}$ of $\eta=\{(1,1),(2,0),(3,1),(4,1)\}$ in $\Tau_5^2$}
\end{center}

%

\begin{rmk}\label{rmk:bottom antichain}
For each $\eta\in\Tau_\alpha^2$, the bottom antichain of $\eta$ always exists and is unique. By finite support, there exist only finitely many $\beta_0,...,\beta_n$ with $\beta_0<\cdots<\beta_n<\alpha$ such that $\eta(\beta_i)=1$ for each $i\le n$. For each $i\le n$, let $\xi_i=\eta|_{[\beta_i+1,\alpha)}\cup\{(\beta_i,0)\}$
Then $\xi_0,...,\xi_n$ form the bottom antichain of $\eta$. Thus the cardinality of the bottom antichain of a given node is always finite. Also note that $\{\nu\in\Tau_\alpha^2:\nu\perp_{lex}\eta\}$ is the union of cones at the nodes in the bottom antichain of $\eta$. If there is no $\beta$ such that $\eta(\beta)=1$, then the bottom antichain of $\eta$ is $\emptyset$.
\end{rmk}

\medskip


\begin{lem}\label{lem:a sequence of pre-antichain tree coheirs}
Let $T$ be any theory, $\M$ a monster model of $T$, $M$ a small model of $T$, $a\in\M$  a tuple of parameters, $p(x)$ a global invariant type over $M$, and $q(x)$ a global coheir over $M$ with $q|_M=p|_M=\tp(a/M)$. Let $\text{ON}$ be the class of all ordinal numbers. Then we can continue constructing a sequence of global coheirs $(q_\alpha((x_\eta)_{\eta\in\Tau_\alpha^2}))_{\alpha\in \text{ON}}$ over $M$ such that:
\begin{itemize}
\item[(i)] $q_0( (x_\eta)_{\eta\in\Tau_0^2})=q((x_\eta)_{\eta\in\Tau_0^2})$.
\item[(ii)] $q_\beta\subseteq q_\alpha$ for all $\beta<\alpha$ with respect to the canonical inclusion.
\item[(iii)] For each $\beta<\alpha$, 
\[\bigotimes_{i<\Omega^{\alpha,2}_\beta}q_\beta ((x_\eta)_{\eta\in C^{\alpha,2}_{\eta(\alpha,\beta,i)}})\subseteq q_\alpha. \]
In other words, $\bigotimes_{i<\Omega^{\alpha,2}_\beta}q_\beta ((x_\eta)_{\eta\in C^{\alpha,2}_{\eta(\alpha,\beta,i)}})=q_\alpha|_{\bigcup_{i<\Omega^{\alpha,2}_\beta}C^{\alpha,2}_{\eta(\alpha,\beta,i)}}$. 
So the type obtained by restricting $q_\alpha$ to the sequence of cones at nodes on the $\beta$-th level is a product of $q_\beta$. 
\item[(iv)] For any set $A$, if $(a_\eta)_{\eta\in\Tau_\alpha^2}\models q_\alpha|_{MA}$, then
\begin{itemize}
\item[($\ast$)] $a_\eta\models p|_{Ma_{\rhd\eta}}$ for each $\eta \in \Tau_\alpha^2$,
\item[($\ast\ast$)] for each $\eta\in\Tau_\alpha^2$, and its bottom antichain $\xi_{0}\lex \cdots \lex\xi_n$ in $\Tau_\alpha^2$, we have
\[\;\;\;\;\;\;\;\;\;\;\;\;\;\;\;\;\;\;\;(a_\nu)_{\nu\in C_\eta}(a_\nu)_{\nu\in C_{\xi_n}}\!...(a_\nu)_{\nu\in C_{\xi_0}}\!\models \left(q_\alpha|_{C_\eta}\!\!\otimes q_\alpha|_{C_{\xi_n}}\!\!\otimes \cdots \otimes q_\alpha|_{C_{\xi_0}}\right)|_{MA}.\]
In particular, for each $\beta<\alpha$ and $\eta\in P^{\alpha,2}_\beta$, we have $a_\eta\models q^\beta|_{M\!Aa_{\perp_{lex}\eta}}$ where $q^\beta = q_\beta|_{\{x_\emptyset\}}$ for $\emptyset\in\Tau_\beta^2$.
\end{itemize}
\end{itemize}
\begin{proof}
We use induction on $\alpha\in\text{ON}$. Let $q_0( (x_\eta)_{\eta\in\Tau_0^2})=q((x_\eta)_{\eta\in\Tau_0^2})$. Suppose that $\alpha>0$ and we have constructed a sequence of global coheirs $(q_\beta((x_\eta)_{\eta\in\Tau_\beta^2}))_{\beta < \alpha}$ satisfying (i), (ii), (iii), and (iv).

First we assume $\alpha$ is a successor ordinal, namely $\alpha=\beta+1$. Choose any $(a^0_\eta)_{\eta\in\Tau_\beta^2}$, $(a^1_\eta)_{\eta\in\Tau_\beta^2}$, and $a^*$ such that 
\begin{align*}
(a^0_\eta)_{\eta\in\Tau_\beta^2}&\models q_\beta |_M,\\
(a^1_\eta)_{\eta\in\Tau_\beta^2}&\models q_\beta |_{M(a^0_\eta)_{\eta\in\Tau_\beta^2}},\;\text{and}\\
a^*&\models\;p\;|_{M(a^1_\eta)_{\eta\in\Tau_\beta^2}(a^0_\eta)_{\eta\in\Tau_\beta^2}}.
\end{align*}
Then \[q_\beta( (x^1_\eta)_{\eta\in\Tau_\beta^2})\otimes q_\beta((x^0_\eta)_{\beta\in\Tau_\beta^2})\supseteq \tp((a^1_\eta)_{\eta\in\Tau_\beta^2}(a^0_\eta)_{\eta\in\Tau_\beta^2}/M).\] Choose any $(\hat{a}^1_\eta)_{\eta\in\Tau_\beta^2}$ and $(\hat{a}^0_\eta)_{\eta\in\Tau_\beta^2}$ outside the monster model $\M$ such that \[(\hat{a}^1_\eta)_{\eta\in\Tau_\beta^2}(\hat{a}^0_\eta)_{\eta\in\Tau_\beta^2}\models q_\beta( (x^1_\eta)_{\eta\in\Tau_\beta^2})\otimes q_\beta((x^0_\eta)_{\beta\in\Tau_\beta^2}).\] Then $(\hat{a}^1_\eta)_{\eta\in\Tau_\beta^2}(\hat{a}^0_\eta)_{\eta\in\Tau_\beta^2}\ind^u_M\M$ by strong left transitivity of $\ind^u$. Since \[(\hat{a}^1_\eta)_{\eta\in\Tau_\beta^2}(\hat{a}^0_\eta)_{\eta\in\Tau_\beta^2}\equiv_M({a}^1_\eta)_{\eta\in\Tau_\beta^2}({a}^0_\eta)_{\eta\in\Tau_\beta^2},\] we have $a^{**}$ such that \[(\hat{a}^1_\eta)_{\eta\in\Tau_\beta^2}(\hat{a}^0_\eta)_{\eta\in\Tau_\beta^2}a^{**}\equiv_M({a}^1_\eta)_{\eta\in\Tau_\beta^2}({a}^0_\eta)_{\eta\in\Tau_\beta^2}a^*.\]
By left extension of $\ind^u$, there exists $\hat{a}^*\equiv_{M(\hat{a}^1_\eta)_{\eta\in\Tau_\beta^2}(\hat{a}^0_\eta)_{\eta\in\Tau_\beta^2}} a^{**}$ such that
\[(\hat{a}^1_\eta)_{\eta\in\Tau_\beta^2}(\hat{a}^0_\eta)_{\eta\in\Tau_\beta^2}\hat{a}^*\ind^u_M\M.\] 
Note that $\hat{a}^*\notin \M$ and $(\hat{a}^1_\eta)_{\eta\in\Tau_\beta^2}(\hat{a}^0_\eta)_{\eta\in\Tau_\beta^2}\hat{a}^*\equiv_M({a}^1_\eta)_{\eta\in\Tau_\beta^2}({a}^0_\eta)_{\eta\in\Tau_\beta^2}a^*.$
For each $\eta\in\Tau_\alpha^2$, define $a^\dagger_\eta$ by
\begin{align*}
a^\dagger_\eta=
\begin{cases}
\hat{a}^* & \text{if } \eta=\emptyset \\
\hat{a}^0_\nu & \text{if } \eta=(0)\coc\nu \\
\hat{a}^1_\nu & \text{if } \eta=(1)\coc\nu, \\
\end{cases}
\end{align*}
and let \[q_\alpha((x_\eta)_{\eta\in\Tau_\alpha^2})=\tp((a^\dagger_\eta)_{\eta\in\Tau_\alpha^2} /\M).\]
Then by Fact \ref{fact:type product} and the induction hypothesis, $q_\alpha$ satisfies (ii) and (iii). Now we show $q_\alpha$ satisfies (iv). Choose any $A$ and suppose $(a_\eta)_{\eta\in\Tau_\alpha^2}\models q_\alpha|_{MA}$. Then \[(a_\eta)_{\eta\in\Tau_\alpha^2}\equiv_M\!  (a^0_\eta)_{\eta\in\Tau_\beta^2}(a^1_\eta)_{\eta\in\Tau_\beta^2}a^*.\]
Thus ($\ast$) is clear. ($\ast\ast$) is by the induction hypothesis and Fact \ref{fact:type product}.

Now we suppose that $\alpha$ is a limit ordinal. Then we just take $q_\alpha:=\bigcup_{\beta<\alpha}p_\beta$ with respect to canonical inclusions. Clearly that $q_\alpha$ satisfies (ii), (iii), and (iv) by finite support.
\end{proof}
\end{lem}

\begin{thm}\label{thm:Kim-dividing->coheir-dividing}
Suppose that $T$ is NATP and let $M$ be a model of $T$. If $\varphi(x,a)$ Kim-divides over $M$, then it coheir-divides over $M$.
\begin{proof}
Suppose $\varphi(x,a)$ Kim-divides over $M$. Then there exist $k\in\omega$ and a global $M$-invariant type $p(y)$ containing $\tp(a/M)$ such that the set $\{\varphi(x,a_i)\}_{i\in\omega}$ is $k$-inconsistent for all $(a_i)_{i\in\omega}\models p^{\otimes\omega}|_M$. To get a contradiction, we assume that there is no global coheir extension $q(y)$ of $\tp(a/M)$ such that $\{\varphi(x,a_i)\}_{i\in\omega}$ is inconsistent for some (any) $(a_i)_{i\in\omega}\models q^{\otimes\omega}|_M$. Recall that the number of all global coheir extensions of $\tp(a/M)$ is bounded (Remark \ref{rmk:the number of inv tp}). Let $\theta$ be the cardinality of all global coheir extensions of $\tp(a/M)$ and $\kappa$ be a cardinal such that $\kappa>\theta^+$.

Choose any global coheir extension $q(y)$ of $\tp(a/M)$. By Lemma \ref{lem:a sequence of pre-antichain tree coheirs}, we can find $(q_\beta( (y_\eta)_{\eta\in\Tau_\beta^2}))_{\beta\le\kappa}$ such that 
\begin{itemize}
\item[(i)] $q_0( (y_\eta)_{\eta\in\Tau_0^2})=q((y_\eta)_{\eta\in\Tau_0^2})$,
\item[(ii)] $\bigotimes_{i<\Omega^{\kappa,2}_\beta}q_\beta ((y_\eta)_{\eta\in C^{\kappa,2}_{\eta(\kappa,\beta,i)}})\subseteq q_\kappa( (y_\eta)_{\eta\in\Tau_\kappa^2})$ for each $\beta<\kappa$.
\item[(iii)] If $(a_\eta)_{\eta\in\Tau_\kappa^2}\models q_\kappa( (y_\eta)_{\eta\in\Tau_\kappa^2})|_M$, then
\begin{itemize}
\item[($\ast$)] $a_\eta\models p|_{Ma_{\rhd\eta}}$ for each $\eta \in \Tau_\kappa^2$,
\item[($\ast\ast$)] for each $\beta<\kappa$, there exists a global coheir $q^{\beta}$ such that $a_\eta\models q^\beta|_{Ma_{\perp_{lex}\eta}}$ for all $\eta\in P^{\kappa,2}_\beta$.
\end{itemize}
\end{itemize}
Let $(a_\eta)_{\eta\in\Tau_\kappa^2}$ be a realization of $q_\kappa|_M$. Since $\kappa>\theta^+$, there exists a set of successor ordinals $K\subseteq\kappa$ such that $|K|=\omega$ and $q^\beta=q^{\beta'}$ for all $\beta,\beta'\in K$. Let $\{\beta_i\}_{i<\omega}$ be an enumeration of $K$ such that $\beta_i<\beta_j$ for all $i<j<\omega$. For each $\eta\in\Tau_\omega^2$, let $f(\eta)$ be an element in $\Tau_\kappa^2$ such that $\text{dom}(f(\eta))=[\beta_{t(\eta)},\kappa)$ and
\begin{align*}
f(\eta)(\beta)=
\begin{cases}
\eta(i) & \text{if } \beta=\beta_{i+1}-1 \\
0 & \text{otherwise}. \\
\end{cases}
\end{align*}
Then we can regard $f$ as an $L_1$-embedding from $\Tau_\omega^2$ into $\Tau_\kappa^2$. Let $q^*=q^\beta$ for some (any) $\beta\in K$. Then $a_{f(\eta)}\models q^*|_{Ma_{\perp_{lex}f(\eta)}}$ for all $\eta\in\Tau_\omega^2$.
 For each $\eta\in\Tau_\omega^2$, let $a^*_\eta:=a_{f(\eta)}$. Then for each antichain $X$ in $\Tau_\omega^2$, $(a^*_\eta)_{\eta\in X}$ is a coheir Morley sequence over $M$ generated by $q^*$ with respect to $\lex$. Thus the set $\{\varphi(x,a^*_\eta)\}_{\eta\in X}$ is consistent. On the other hand, if $X$ is a path in $\Tau_\omega^2$, then $(a^*_\eta)_{\eta\in X}$ is an invariant Morley sequence over $M$ generated by $p$ with respect to $\lhd$. Thus the set $\{\varphi(x,a^*_\eta)\}_{\eta\in X}$ is $k$-inconsistent. By compactness and \cite[Lemma 3.20]{AKL21}, we can construct an antichain tree, which yields a contradiction with the assumption that $T$ is NATP.
\end{proof}
\end{thm}

\begin{rmk}
More generally, by the same argument in Lemma \ref{lem:a sequence of pre-antichain tree coheirs} and Theorem \ref{thm:Kim-dividing->coheir-dividing}, we can say that if there exists a pre-independence relation $\ind$ stronger than $\ind^i$, satisfying monotonicity, strong finite character, strong left transitivity, right extension, and left existence over a model $M$, then for each formula $\varphi(x,a)$ Kim-dividing over $M$, there exists a global invariant type $p(x)\supseteq\tp(a/M)$ such that: 
\begin{itemize}
\item[(i)] $a'\ind_MB$ whenever $a'\models p|_{MB}$,
\item[(ii)] $\{\varphi(x,a_i)\}_{i<\omega}$ is inconsistent for all $(a_i)\models p^{\otimes\omega}|_M$.
\end{itemize}
\end{rmk}

From Theorem \ref{thm:Kim-dividing->coheir-dividing}, we get some observations in NATP theories, that correspond to some phenomena in NTP$_2$ theories appeared in \cite{CK12}.
First we recall a special case of the Broom Lemma \cite[Lemma 3.1]{CK12}.

\begin{fact}\label{fact:Broom Lemma}
Let $M$ be a model. Suppose that $\alpha(x,e)\vdash \psi(x,c)\vee \bigvee_{i<n}\varphi_i(x,a_i)$ and $\varphi_i(x,a_i)$ coheir-divides over $M$ for each $i<n$. Then there exists $e_0,...,e_m$ such that $e_i\equiv_M e$ for each $i\le m$ and $\{\alpha(x,e_i)\}_{i\le m}\vdash \psi(x,c)$.
\end{fact}

As forking implies quasi-dividing in NTP$_2$ theories, Kim-forking implies quasi-dividing over models in NATP theories.

\begin{cor}\label{cor:Kim-forking >> quasi-divide}
Suppose that $T$ is NATP. If a formula Kim-forks over a model $M$, then it quasi-divides over $M$.
\begin{proof}
Suppose a formula $\varphi(x,a)$ Kim-forks over a model $M$. Then there exist formulas $\psi_0(x,a_0), ... ,\psi_n(x,a_n)$ such that $\varphi(x,a)\vdash\bigvee_{i\le n}\psi_i(x,a_i)$ and $\psi_i(x,a_i)$ Kim-divides over $M$ for each $i\leq n$. By Lemma \ref{thm:Kim-dividing->coheir-dividing}, $\psi_i(x,a_i)$ coheir-divides over $M$ for each $i\le n$. By taking $\psi(x,c):=\perp$ ({\it i.e.}, $\forall x(x\neq x)$), we can apply Fact \ref{fact:Broom Lemma}. Thus there exist $a_0,...,a_m$ such that $a_i\equiv_M a$ for each $i\le m$ and $\{\varphi(x,a_i)\}_{i\le m}\vdash \perp$. So $\{\varphi(x,a_i)\}_{i\le m}$ is inconsistent and $\varphi(x,a)$ quasi-divides over $M$.
\end{proof}
\end{cor}

The strategy of the proof of Corollary \ref{cor:pre_existence of strong coheir} is from \cite[Proposition 3.7]{CK12}. To generalize the statement, we need one more property for pre-independence relations.

\begin{dfn}\cite{Mut22}\label{def:property_of_ind2}
 We say a pre-independence relation $\ind$ satisfies {\it quasi-strong finite character (over $d$)} if for each $a,b$, there exists a partial type $\Sigma(x,y)$ such that $a'b'\models \Sigma(x,y)$ if and only if $a'\equiv_da$, $b'\equiv_db$, $a'\ind_db'$.
\end{dfn}

Suppose $\ind$ satisfies quasi-strong finite character over $d$. Note that for each $a,b$, if $\{\Sigma_i(x,y)\}_{i\in I}$ is a set of partial types satisfying the condition in Definition \ref{def:property_of_ind2}, then $\bigcup_{i\in I}\Sigma_i(x,y)$ also satisfies the condition. And if $\Sigma(x,y)$ satisfies the condition, then it is $d$-invariant hence it is $d$-definable. So we can consider the maximality of such types definable over $d$. 

\begin{nota}
Suppose that $\ind$ satisfies quasi-strong finite character and full existence over $d$. For each $a$ and $b$, let $\Sigma_{a\inds_db}(x,y)$ be the unique $\subseteq$-maximal partial type over $d$ such that \[a'b'\models \Sigma_{a\inds_db}(x,y)\] if and only if \[a'\equiv_da,\; b'\equiv_db, \text{ and } a'\ind_db'.\]

Note that if $\ind$ satisfies monotonicity over $d$ additionally, then $\Sigma_{a'\inds_db'}(x',y')\subseteq\Sigma_{a\inds_db}(x,y)$ for all $a'\subseteq a$ and $b'\subseteq b$, by the maximality.

\end{nota}


\begin{cor}\label{cor:pre_existence of strong coheir}
Suppose that $T$ is NATP and let $\ind$ be a pre-independence relation satisfying monotonicity, quasi-strong finite character, and full existence over a model $M$. Let $a\in\M$ be a tuple of parameters. Then there exists a global type $p(x)$ containing $\tp(a/M)$ such that $a'\ind_M A$ and $A\ind_M^K a'$ for all $a'\models p(x)|_{MA}$.
\begin{proof}
Let $q(x):=\tp(a/M)$ and $\Lambda(x)$ be
\begin{equation*}
\begin{aligned}
q(x) {} & \cup \{ \neg\varphi(x',b')\;|\;\varphi(x',y')\in\L(M),\;b'\in\M, \;x'\text{ is a finite subtuple of }x,\\ & \qquad\qquad\qquad\;\; \varphi(a',y')\text{ Kim-forks over }M \text{ for the finite subtuple }a'\subseteq a\\ 
& \qquad\qquad\qquad\;\;\text{corresponding to }x'\subseteq x\} \\
      & \cup \{ \neg\psi(x''\!,d'')\,|\,\psi(x''\!,z'')\in\L(M),\;d''\!\!\in\!\M,\; x''\text{ is a finite subtuple of }x,\\
       &\qquad\qquad\qquad\;\;\neg\psi(x'',z'')\in\Sigma_{a\inds_Md}(x,z)\text{ for some }d\supseteq d''\text{ and }z\supseteq z''\} 
       \\
       &\qquad\qquad\qquad\;\;\text{where } d''\text{ is a subtuple of }d \text{ corresponding to } z''\subseteq z\}
\end{aligned}
\end{equation*}
First we show that $\Lambda(x)$ is consistent. Suppose not. Then there exist $\varphi_0(x'_0,b'_0)$, $...$ , $\varphi_n(x'_n,b'_n)$ and $\psi_0(x''_0,d''_0)$, $...$ , $\psi_m(x''_m,d''_m)$ such that $\varphi_i(a'_i,y'_i)$ Kim-forks over $M$ and $a_i'\subseteq a$ for each $i\le n$, $\neg\psi_j(x''_j,z''_j)\in\Sigma_{a\inds_Md_j}(x,z_j)$ for each $j\le m$, and 
\[
q(x)\vdash \bigvee_{i\le n}\varphi_i(x'_i,b'_i) \vee \bigvee_{j\le m}\psi_j(x''_j,d''_j).\]
Note that $x'_0,...,x'_n$, $x''_0,...,x''_m$ are finite subtuples of $x$, and $d''_j\subseteq d_j$ for each $j\le m$. Since $\bigvee_{i\le n}\varphi_i(a'_i,y'_i)$ also Kim-forks over $M$, we may assume $n=0$. Let  $a':=a'_0$, $b':=b'_0$, $x':=x_0'$, $y':=y_0'$, and $\varphi(x',y'):=\varphi_0(x'_0,y'_0)$. Then we have 
\[
q(x)\vdash \varphi(x',b') \vee \bigvee_{j\le m}\psi_j(x''_j,d''_j).\]

By Corollary \ref{cor:Kim-forking >> quasi-divide}, $\varphi(a',y')$ quasi-divides over $M$, so there exist $a_0,...,a_k$ such that $\{\varphi(a_i,y)\}_{i\le k}$ is inconsistent and $a_i\equiv_M a'$ for each $i\le k$. Choose any $\hat{a}_i$ containing $a_i$ such that $\hat{a}_ia_i\equiv_Maa'$ for each $i\le k$. Let $\hat{a}:=(\hat{a}_0,...,\hat{a}_k)$ and $r(\hat{x}_0\cdots\hat{x}_k):=\tp(\hat{a}_0\cdots\hat{a}_k/M)$. For each $i\le k$ and $j\le m$, let $\hat{x}'_i$ and $\hat{x}''_{ij}$ be the subtuples of $\hat{x}_i$ whose indices in $\hat{x}_i$ correspond to the indices of $x'$ and $x''_j$ in $x$, respectively.
Then for each $i\le k$,
\[ r|_{\hat{x}_i}\vdash \varphi(\hat{x}'_i,b)\vee \bigvee_{j\le m}\psi_j(\hat{x}''_{ij},d''_j).\]
Thus 
\[ r(\hat{x}_0,...,\hat{x}_k)\vdash\bigwedge_{i\le k}\Big[ \varphi(\hat{x}'_i,b)\vee \bigvee_{j\le m}\psi_j(\hat{x}''_{ij},d''_j)\Big].\]
But 
\[r(\hat{x}_0,...,\hat{x}_k)\vdash \neg\exists y\Big(\bigwedge_{i\le k}\varphi(\hat{x}'_i,y)\Big),\]
and hence
\[r(\hat{x}_0,...,\hat{x}_k)\vdash \bigvee_{i\le k, j\le m}\!\!\!\!\psi_j(\hat{x}''_{ij},d''_j).\] 

By full existence of $\ind$, there exists $a^*=(a^*_0,...,a^*_k)$ such that $a^*\equiv_M\hat{a}$ and $a^*\ind_M d_{\le m}$. Since $a^*\models r$, there exist $i\le k$ and $j\le m$ such that $\models\psi_j(a^*_{ij},d''_j)$, where $a^*_{ij}$ is the subtuple of $a^*_i$ whose indices in $a^*_i$ correspond to the indices of $\hat{x}''_{ij}$ in $\hat{x}_i$ (equivalently, $x''_{j}$ in $x$). Thus $(a^*_i,d_j)\not\models\Sigma_{a\inds_Md_j}(x,z_j)$. Since $a^*_i\equiv_M \hat{a}_i\equiv_M a$ and $d_j\equiv_M d_j$, we have $a^*_i\not\ind_M d_j$, which yields a contradiction with $a^*\ind_M d_{j\le m}$ by monotonicity of $\ind$. Thus $\Lambda(x)$ is consistent.

 Choose any global type  $p(x)$ which is a completion of $\Lambda(x)$. Then $p(x)$ satisfies all the conditions we want.
\end{proof}
\end{cor}

\begin{cor}\label{cor:existence of strong coheir}
Suppose the theory is NATP. Then for each model $M$ and a tuple of parameters $b$, there exists a global coheir $p(x)$ over $M$ containing $\tp(b/M)$ such that $B\ind^K_M b'$ for all $b'\models  p(x)|_{MB}$.
\begin{proof}
Apply Corollary \ref{cor:pre_existence of strong coheir} on $\ind^u$.
\end{proof}
\end{cor}

\section{Some remarks on witnesses of Kim-dividing in NATP theories}

In this section, we discuss the concept of witnesses of Kim-dividing, whose existence in NSOP$_1$ theories is given by Kim's lemma for Kim-dividing in \cite{KR20}. Note that the existence of witnesses of Kim-dividing over models in NTP$_2$ theories is a special case of \cite[Theorem 3.2]{CK12}.

\begin{dfn}\cite[Definition 7.8]{KR20}\label{def:witness}
Let $b$ be a tuple of parameters and $C$ be a set. 
\begin{itemize}
\item[(i)] An indiscernible sequence $(b_i)_{i<\omega}$ over $C$ is called a {\it witness of dividing of $b$ over $C$} if $b_0=b$ and $\{\varphi(x,b_i)\}_{i<\omega}$ is inconsistent whenever $\varphi(x,b)$ divides over $C$.
\item[(ii)] An indiscernible sequence $(b_i)_{i<\omega}$ over $C$ is called a {\it witness of Kim-dividing of $b$ over $C$} if $b_0=b$ and $\{\varphi(x,b_i)\}_{i<\omega}$ is inconsistent whenever $\varphi(x,b)$ Kim-divides over $C$.
\item[(iii)] A global invariant type $p(x)\supseteq\tp(b/C)$ is called a {\it witness of dividing of $b$ over $C$} if every invariant Morley sequence $(b_i)_{i<\omega}$ generated by $p(x)$ over $C$ ({\it i.e.}, $b_i\models p(x)|_{Cb_{<i}}$ for each $i<\omega$) is a witness of dividing of $b$ over $C$.
\item[(iv)] A global invariant type $p(x)\supseteq\tp(b/C)$ is called a {\it witness of Kim-dividing of $b$ over $C$} if every invariant Morley sequence $(b_i)_{i<\omega}$ generated by $p(x)$ over $C$ is a witness of Kim-dividing of $b$ over $C$.
\end{itemize}
If $a$ and $B$ are clear in the context, then we just say the sequence (type) is a witness of dividing or Kim-dividing.
\end{dfn}

\begin{dfn}
A global type $p(x)$ is said to be {\it strictly invariant with respect to forking over $C$} (or we just say that the type is {\it strictly invariant over $C$}) if $a\ind^{
i}_CB$ and $B\ind^{f}_Ca$ whenever $a\models p|_{CB}$. We write $a\ind^{ist}_Cb$ if there exists a strictly invariant global type $p(x)$ containing $\tp(a/Cb)$. A sequence $(b_i)_{i<\omega}$ is said to be {\it strictly invariant over $C$} if it is generated by a strictly invariant global type over $C$.
\end{dfn}

The following statements are consequences of Kim's lemma for dividing or Kim-dividing, which appear in \cite{Kim98}, \cite{KR20}, and \cite{CK12} 

\begin{fact}\cite{Kim98}\cite{KR20}\cite{CK12}
Let $b$ be a tuple of parameters, $C$ a set, and $M$ a model.
\begin{itemize}
\item[(i)] If $T$ is simple, then every Morley sequence in $\tp(b/C)$ is a witness of dividing of $b$ over $C$
\item[(ii)] If $T$ is NSOP$_1$, then every invariant Morley sequence over $M$ starting with $b$ is a witness of Kim-dividing of $b$ over $M$.
\item[(iii)] If $T$ is NTP$_2$, then every strictly invariant Morley sequence over $M$ starting with $b$ is a witness of dividing of $b$ over $M$.
\end{itemize}
\end{fact}

From now we investigate the possibility of the existence of witnesses of Kim-dividing in NATP theories. First we note that the existence of witnesses of Kim-dividing gives a necessary condition of being a witness of Kim-dividing for coheirs.

\begin{dfn}
A global type $p(x)$ is called a {\it strict coheir with respect to Kim-forking over $B$} (or we just call the type a {\it strict coheir over $C$}) if $a\ind^{u}_BC$ and $C\ind^{K}_Ba$ whenever $a\models p|_{BC}$. We write $a\ind^{ust}_Cb$ if there exists a strict coheir $p(x)$ containing $\tp(a/Cb)$. A sequence $(b_i)_{i<\omega}$ is called a {\it strict coheir Morley sequence over $C$} if it is generated by a strict coheir over $C$.
\end{dfn}

\begin{rmk}
By Corollary \ref{cor:existence of strong coheir}, every model is an extension base for $\ind^{ust}$ in NATP theories. 
\end{rmk}

\begin{lemma}\label{lem:equivalent conditions for uKst-ind}
The following are equivalent.
\begin{itemize}
\item[(i)] $a\ind_D^{ust}b$.
\item[(ii)] For all $c$, there exists $c'\equiv_{Db}c$ such that $bc'\ind_D^{K}a$ and $a\ind_D^{u}bc'$.
\end{itemize}
\begin{proof}
Clearly (i) implies (ii). Suppose (ii). Then by compactness, the partial type
\begin{equation*}
\begin{aligned}
\Lambda(x):=\tp_x(a/Db) {} & \cup \{ \neg\varphi(x,c)\;|\;c\in\M,\;\varphi\in\L(D),\;\varphi(a,y)\text{ Kim-forks over }D\} \\
      & \cup \{ \neg\psi(x,c)\;|\;c\in\M,\;\psi\in\L(D),\;\psi(x,c)\text{ is not realized in }D\}
\end{aligned}
\end{equation*}
is consistent. Any global type which is a completion of $\Lambda(x)$ is a strict coheir containing $\tp(a/Db)$.
\end{proof}
\end{lemma}

\begin{prop}\label{prop:necessary cond}
Let $T$ be an NATP theory and $\M$ its monster model. Suppose that for any model $M$ and $a\in\M$, there exists an invariant Morley sequence which is a witness of Kim-dividing of $a$ over $M$. Let $M$ be a model, $a\in\M$, and $I=(a_i)_{i<\omega}$ a coheir Morley sequence over $M$ with $a_0=a$. If $I$ is a witness of Kim-dividing of $a$ over $M$, then it is a strict coheir Morley sequence of $\tp(a/M)$.
\begin{proof}
Let $p$ be a global coheir of $\tp(a/M)$ such that $I\models p^{\otimes\omega}|_M$. Note that the number of all coheirs of $\tp(a/M)$, say $\kappa$, is bounded. Choose any $\kappa'>\kappa$. Let $I'=(a'_i)_{i<\kappa'}$ be a coheir Morley sequence generated by $p$ over $M$. Then $I'$ is also a witness of Kim-dividing of $a$ over $M$. First we claim that $a'_i\ind_M^{ust}a'_{<i}$ for each $i<\kappa'$. Choose any $i<\kappa'$. By Lemma \ref{lem:equivalent conditions for uKst-ind}, it is enough to show that for any $b$, there exists $b'\equiv_{Ma'_{<i}}b$ such that $a'_{<i}b'\ind_M^{K}a'_i$ and $a'_i\ind_M^{u} a'_{<i}b'$. Let $\Theta=\kappa'\setminus i$. Then we can find $J=(a^*_j)_{j\in\Theta}$ such that $a^*_j\models p|_{Ma'_{<i}ba^*_{<j}}$ for all $j\in\Theta$. Since $J\equiv_{Ma'_{<i}}(a'_j)_{j\geq i}$, we can find $b'\equiv_{Ma'_{<i}}b$ such that $a'_j\models p|_{Ma'_{<j}b'}$ for all $i\leq j<\kappa'$. Then $b'$ satisfies $a'_i\ind_M^{u} a'_{<i}b'$. Note that $(a'_j)_{j\in\Theta}$ is a witness of Kim-dividing of $a'_i$ over $M$.
Thus $a'_{<i}b'\ind_M^{Kd}a'_i$ by indiscernibility of $(a'_j)_{j\in\Theta}$ over $a'_{<i}b'$. Note that Kim-dividing and Kim-forking over a model are equivalent since we assume that every type over a model has a global invariant extension that is a witness. Thus $a'_{<i}b'\ind_M^{K}a'_i$.

 Since $\kappa'$ is sufficiently large, we may assume that there exists a global strong coheir $q$ of $\tp(a/M)$ such that $\tp(a_i/Ma_{<i})\subseteq q$ for all $i\in\kappa'$. Thus $I$ is a strict coheir Morley sequence in $\tp(a/M)$.
\end{proof}
\end{prop}

Now we consider a sufficient condition for the existence of witnesses of Kim-dividing. The following statement is analogous to \cite[Theorem 3.11]{CK12}.

\begin{prop}\label{prop:sufficient cond}
Let $T$ be an NATP theory and $M$ a model. Then (i) implies (ii) where:
\begin{itemize}
\item[(i)] There exists a pre-independence relation $\ind$ which is stronger than $\ind^{h}$ and satisfies monotonicity, full existence, the right chain condition for coheir Morley sequences, and strong right transitivity over $M$.
\item[(ii)] For all $b$, there exists a global coheir which is a witness of Kim-dividing of $b$ over $M$. Thus Kim-forking is equivalent with Kim-dividing over $M$.
\end{itemize}
\begin{proof}
Suppose (i) is true and choose any $b$. By full existence, there exists a global type $p(x)\supseteq\tp(b/M)$ such that $A\ind_Mb'$ for all $b'\models p|_{MA}$. Since $\ind$ is stronger than $\ind^h$, $p(x)$ is a coheir over $M$. Let $\varphi(x,b)$ Kim-divide over $M$. Then by Theorem \ref{thm:Kim-dividing->coheir-dividing}, there exists a global coheir $q(x)$ over $M$ and $k<\omega$ such that $\{\varphi(x,b_i)\}_{i<\omega}$ is $k$-inconsistent for every coheir Morley sequence $(b_i)_{i<\omega}$ generated by $q(x)$ over $M$. 

It is enough to show that $\{\varphi(x,b_i)\}_{i<\omega}$ is inconsistent for any coheir Morley sequence $(b_i)_{i<\omega}$ generated by $p(x)$ over $M$. 
To get a contradiction, suppose not. Then $\{\varphi(x,b_i)\}_{i<\omega}$ is consistent for some (any) coheir Morley sequence $(b_i)_{i<\omega}$ generated by $p(x)$ over $M$.
 \smallskip
\begin{claim}

For each $n<\omega$ and a small set $A$, there exists $(b_\eta)_{\eta\in 2^{<n}}$ such that:
\begin{itemize}
\item[($\ast$)] $b_\eta \models q|_{Mb_{\rhd\eta}}$ for each $\eta\in 2^{<n}$,
\item[($\ast\ast$)] $b_\eta \models p|_{MAb_{\perp_{lex}\eta}}$ for each $\eta\in 2^{<n}$, and 
\item[($\ast\!\ast\!\ast$)] $A\ind_M \!(b_\eta)_{\eta\in2^{<n}}$.
\end{itemize}
\smallskip
\noindent{\it Proof of Claim.} We use induction on $n<\omega$. If $n=0$, then there is nothing to prove. Suppose $n=1$ and let $A$ be an arbitrary small set. Choose any $b'\models p|_{MA}$ and let $b_{\emptyset}:=b'$. then $(b_\eta)_{\eta\in2^{<n}}$ satisfies ($\ast$), ($\ast\ast$), and ($\ast\!\ast\!\ast$).

Now let $n$ be an arbitrary natural number larger than $0$ and suppose that for any small set, there exists $(b_i)_{\eta\in2^{<n}}$ satisfying ($\ast$), ($\ast\ast$), and ($\ast\!\ast\!\ast$) over the set. Choose any $A$. We find $(b_i)_{\eta\in2^{<{n+1}}}$ satisfying ($\ast$), ($\ast\ast$), and ($\ast\!\ast\!\ast$) over $A$.

By the induction hypothesis, there exists $(b^0_\eta)_{\eta\in2^{<n}}$ satisfying ($\ast$), ($\ast\ast$), and ($\ast\!\ast\!\ast$) over $A$. Let $B^0:=(b^0_\eta)_{\eta\in 2^{<n}}$. By applying the induction hypothesis again, we can find $(b^1_i)_{\eta\in2^{<n}}$ satisfying ($\ast$), ($\ast\ast$), and ($\ast\!\ast\!\ast$) over $AB^0$. Let $B^1:=(b^1_\eta)_{\eta\in 2^{<n}}$. Then $A\ind_M B^0B^1$ by strong right transitivity. 

Let $m=|2^{<{n+1}}|-1=2^{n+1}-2$ and choose any enumeration $\{\eta_i\}_{i<m}$ of $2^{<{n+1}}\setminus\{\emptyset\}$ with $\eta_0(i)=0$ for all $i<n+1$. Let $B=(b_{\eta_i})_{i<m}$ be an enumeration of $B^0B^1$ such that $b_{\eta_0}=b^0_{\eta^*}$ where $\eta^*\in2^{n-1}$ and $\eta^*(i)=0$ for all $i<n-1$. Then $A\ind_MB$ and $b_{\eta_0}\models q|_M$.
\medskip
\begin{subclaim}
We can continue constructing a sequence
\[B_0,B_1,...,B_\alpha,B_{\alpha+1},...\]
such that $B_\alpha=(b_{\alpha,\eta_0},...,b_{\alpha,\eta_{m-1}})$ for each ordinal $\alpha$, $b_{0,\eta_i}=b_{\eta_i}$ for each $i<m$ so that $B_0=B$, and
\begin{itemize}
\item[($\dagger$)] $B_\alpha\equiv_MB_0$ and $B_\alpha\ind_M^uB_{<\alpha}$ for each ordinal $\alpha$,
\item[($\ddagger$)] $b_{\alpha,\eta_0}\models q|_{MB_{<\alpha}}$ for each ordinal $\alpha$.
\end{itemize}
\noindent{\it Proof of Subclaim.}
Let $b_{0,\eta_i}:=b_{\eta_i}$ for each $i<m$ and $B_0:=(b_{0,\eta_0},...,b_{0,\eta_{m-1}})$. For each ordinal $\alpha$, suppose that we have constructed a sequence $(B_\beta)_{\beta<\alpha}$ satisfying ($\dagger$) and ($\ddagger$). Choose any $b'\models q|_{MB_{<\alpha}}$. Since $b'\equiv_M b_{0,f\eta_0}$, there exists $B'$ containing $b'$ such that $B'\equiv_MB_0$. Since $b'\ind^u_M B_{<\alpha}$, there exists $B''\equiv_{Mb'}B'$ such that $B''\ind^u_M B_{<\alpha}$ by left extension of $\ind^u$. Let $B_\alpha=B''$ and $b_{\alpha,\eta_0}=b'$. Then the sequence $(B_\beta)_{\beta<\alpha+1}$ satisfies ($\dagger$) and ($\ddagger$). This completes proof of the subclaim.
\qed
\end{subclaim}
\smallskip
Let $\kappa$ be a sufficiently large cardinal and let $(B_\alpha)_{\alpha<\kappa}$ be a sequence given by the subclaim. Since the number of global coheirs over $M$ containing $\tp(B/M)$ is bounded, we may assume that the sequence is a coheir Morley sequence over $M$ with $B_0=B$. By right chain condition for coheir Morley sequences, we may assume $A\ind_M(B_\alpha)_{\alpha<\kappa}$ and $(B_\alpha)_{\alpha<\kappa}$ is indiscernible over $MA$. By monotonicity of $\ind$, we have $A\ind_M Bb_{1,\eta_0}$. By the indiscernibility, we have $B_1\equiv_{MA}B_0=B$, in particular $b_{1,\eta_0}\equiv_{MA}b_{0,\eta_0}=b_{\eta_0}$. Thus $b_{1,\eta_0}\models p|_{MA}$. Let $b_\emptyset:=b_{1,\eta_0}$. Then $(b_\eta)_{\eta\in 2^{<n+1}}$ satisfies ($\ast$), ($\ast\ast$), and ($\ast\!\ast\!\ast$). This completes proof of the claim.\qed
\end{claim}
\smallskip
 As in Theorem \ref{thm:Kim-dividing->coheir-dividing}, we can construct a $k$-antichain tree by compactness. Thus $T$ is NATP by \cite[Lemma 3.20]{AKL21}, it is a contradiction.
\end{proof}
\end{prop}

\begin{rmk}\label{rmk:ck-indep mut22}
In \cite[Proposition 3.1]{Mut22}, Mutchnik shows that in any theory, there exists a pre-independence relation which is called $\ind^{\rm CK}$, stronger than $\ind^h$, and satisfies monotonicity, full existence, and right chain condition for coheir Morley sequences over models. Moreover it satisfies right extension and quasi-strong finite character over models. 

But $\ind^{\rm CK}$ does not satisfy strong right transitivity in general. Let $M$ be a small model of DLO, $a,b,c$ live in the same cut of $M$, and $a<b<c$. Then $b\ind^{\rm CK}_M a$, $ba\ind^{\rm CK}_Mc$, but $b\not\ind^{\rm CK}_M ac$. We will explain this in Remark \ref{rmk:proof of remark 4.9}.

\end{rmk}

\section{Some remarks on N-$\omega$-DCTP$_2$}

As we mentioned above, in his work Mutchnik proved that SOP$_1$ and SOP$_2$ are equivalent at the level of theories. This is a surprising result in itself, but the technique he used to prove it is also interesting. We end this paper by leaving some remarks on his works. 

\begin{dfn}\cite{Mut22}
Let $\ind$ be a pre-independence relation satisfying monotonicity, right extension, quasi-strong finite character, and full extension over models. Let $M$ be a model and $p(x)\in S(M)$. We say a formula $\varphi(x,b)$ {\it $h^{\inds}$-coheir-divides with respect to a type $p(x)$} if there is a coheir Morley sequence $(b_i)_{i<\omega}$ over $M$ with $b_0=b$ such that no $a$ satisfies $a \ind_M (b_i)_{i<\omega}$ and $a\models p(x)\cup\{\varphi(x,b_i)\}_{i<\omega}$. We say a formula $\varphi(x,b)$ {\it $h^{\inds}$-coheir-forks with respect to a type $p(x)$} if there are $\psi_0(x,b_0),...,\psi_n(x,b_n)$ such that $\varphi(x,b)\vdash \bigvee_{i\le n}\psi_i(x,b_i)$ and $\psi_i(x,b_i)$ $h^{\inds}$-coheir-divides with respect to $p(x)$ for each $i\le n$.  For a given $\ind$, define a pre-independence relation $\ind'$ as follows: $a\ind'_Mb$ if and only if $\varphi(x)$ does not $h^{\inds}$-coheir fork with respect to $\tp(a/M)$ for all $\varphi(x)\in\tp(a/Mb)$

Let $\ind^0:=\ind^h$, and $\ind^{n+1}:=(\ind^n)'$ for each $n<\omega$. Let $\ind^{\rm CK}:=\bigcap_{n<\omega}\ind^n$.
\end{dfn}

\begin{rmk}\label{rmk:proof of remark 4.9}
We prove the second paragraph of Remark \ref{rmk:ck-indep mut22}. $b\not\ind^{\rm CK}_M ac$ is clear since $a<x<c$ is realized by $b$ and it coheir-divides over $M$. 

By using induction on $n<\omega$, we show $D\ind^n_M E$ for all $n<\omega$ and two non-empty small sets $D$ and $E$ such that:
\begin{itemize}
\item[(i)] $d<e$ and $d\equiv_Me$ for all $d\in D$, $e\in E$,
\item[(ii)] there exist $m,m'\in M$ and $d\in D$ such that $m<d<m'$.
\end{itemize} 
If $n=0$, then the case is clear. Now let $n<\omega$ and assume that the statement is true for $n$. Suppose $D$ and $E$ satisfy (i) and (ii) above, and $D\not\ind^{n+1}_ME$. Note that every element in $D\cup E$ lives in the same cut. There exists a formula $\varphi(\bar{x},\bar{e})\in\tp(D/ME)$ such that $\varphi(\bar{x},\bar{e})$ $h^{\inds^{n}}$-coheir-forks with respect to $\tp(D/M)$. We may assume that $\varphi$ is of the form \[m<x_0<\cdots<x_k<e\]
for some $m\in M$ and $e\in E$. By \cite[Corollary 3.7.1]{Mut22}, there exist $e_0,...,e_l$ such that $e_i\equiv_M e$ for each $i\le l$ and there is no $d_0...d_k\models\tp(D/M)$ such that \[d_0... d_k\ind^n_Me_0... e_l\] and \[m<d_0<\cdots<d_k<e_i\] for all $i\le l$. But we can choose such $d_0,...,d_k$ by the induction hypothesis. It is a contradiction.

Thus $ba\ind^{\rm CK}_Mc$. We can show $b\ind^{\rm CK}_M a$ by the same argument. \qed
\end{rmk}

As we mentioned in Remark \ref{rmk:ck-indep mut22}, it is proved that $\ind^{\rm CK}$ is stronger than $\ind^h$ and satisfies monotonicity, right extension, quasi-strong finite character, full existence, and right chain condition for coheir Morley sequences over models \cite{Mut22}.

By full existence, we can define a class of global types stronger than the class of coheirs as follows:

\begin{dfn}\cite{Mut22}
For a model $M$, we say a global type $p(x)$ is a {\it canonical coheir over $M$} if $B\ind^{\rm CK}_M a$ whenever $a\models p(x)|_{MB}$. 

For any given $b$, there exists $b'\equiv_Mb$ such that $\M\ind^{\rm CK}_M b'$ by full existence. Thus $\tp(b'/\M)$ is a global canonical coheir over $M$ containing $\tp(b/M)$.
\end{dfn}

Then by Theorem \ref{thm:Kim-dividing->coheir-dividing}, we have the following observation.

\begin{rmk}
Every canonical coheir is a strict coheir in NATP theories since Kim-dividing implies coheir-dividing in NATP, and coheir-dividing implies $h^{\inds}$-coheir-dividing in any theory, for any pre-independence relation $\ind$.
\end{rmk}

The following tree properties form new dividing lines located inside NATP.

\begin{dfn}\cite{Mut22}\cite{Mut23}\cite{Sim}\label{def:DCTP2}
We say an antichain $X=\{\eta_0,...,\eta_n\}\subseteq \tree$ with $\eta_0\lex \cdots\lex \eta_n$ is a {\it descending comb} if $\eta_i\wedge\eta_k=\eta_j\wedge\eta_k$ for all $i<j<k$.

For $k<\omega$, we say a formula $\varphi(x,y)$ has the {\it k-descending comb tree property 2 (k-DCTP$_2$)} if there is a tree-indexed set $(a_\eta)_{\eta\in\tree}$ of parameters such that
\begin{enumerate}
\item[(i)] $\{\varphi(x,a_\eta)\}_{\eta\in X}$ is consistent for each descending comb $X \subseteq\tree$,
\item[(ii)] $\{\varphi(x,a_{\eta|_i})\}_{i<\omega}$ is $k$-inconsistent for each $\eta\in 2^\omega$.
\end{enumerate}

We say a theory has $k$-DCTP$_2$ if there exists a formula having $k$-DCTP$_2$. We say a theory has $\omega$-DCTP$_2$ if it has $k$-DCTP$_2$ for some $k<\omega$. We say a theory is N-$k$-NDCTP$_2$ (N-$\omega$-DCTP$_2$) if it does not have $k$-DCTP$_2$ ($\omega$-DCTP$_2$).
\end{dfn}

It is easy to see that $k$-DCTP$_2$ is always observed by a strongly indiscernible tree as below.

\begin{rmk}\label{rmk:str-indisc k-DCTP_2}
If a formula $\varphi(x,y)$ has $k$-DCTP$_2$, then there exists a strongly indiscernible tree $(a_\eta)_{\eta\in\tree}$ such that 
\begin{itemize}
\item[(i)] for any descending comb $X\subseteq\tree$, the set $\{\varphi(x,a_\eta)\}_{\eta\in X}$ is consistent,
\item[(ii)] for any $\eta_0,...,\eta_{k-1}$ with $\eta_0\trn ...\trn \eta_{k-1}$, the set $\{\varphi(x,a_{\eta_0}),...,\varphi(x,a_{\eta_{k-1}})\}$ is inconsistent.
\end{itemize}
\begin{proof}
Suppose that a formula $\varphi(x,y)$ has $k$-DCTP$_2$ with $(b_\eta)_{\eta\in\tree}$.
Define a map $f:\treo\to\tree$ by
\[
    f(\eta)= 
\begin{cases}
    \emptyset & \text{ if } \eta=\emptyset\\
    f(\nu)\coc(0)^i\coc (1)       & \text{ if } \eta=\nu\coc(i),
\end{cases}
\]
and let $c_\eta:=b_{f(\eta)}$ for each $\eta\in\treo$. 
By using the modeling property, we have a strongly indiscernible $(d_\eta)_{\eta\in\omega^{<\omega}}$ which is strongly locally based on $(c_\eta)_{\eta\in\omega^{<\omega}}$. Define a map $g:\tree\to\omega^{<\omega}$ by 
\[
    g(\eta)= 
\begin{cases}
    \emptyset & \text{ if } \eta=\emptyset\\
    g(\nu)\coc(0)       & \text{ if } \eta=\nu\coc(1)\\
    g(\nu)\coc(1)       & \text{ if } \eta=\nu\coc(0),
\end{cases}
\]
and let $a_\eta:=d_{g(\eta)}$ for each $\eta\in\tree$. Then $(a_\eta)_{\eta\in\tree}$ is strongly indiscernible and satisfies (i) and (ii).
\end{proof}
\end{rmk}

\begin{rmk} Let $k<\omega$.
\begin{itemize}
\item[(i)] If a theory has $k$-DCTP$_2$, then the theory has TP$_2$.
\item[(ii)] If a theory has $k$-DCTP$_2$, then the theory has SOP$_1$.
\item[(iii)] If a theory has ATP, then the theory has $k$-DCTP$_2$.
\item[(iv)] If a theory has $k$-DCTP$_2$, then the theory has $(k+1)$-DCTP$_2$.
\end{itemize}
\begin{proof}
(iii), (iv) are clear. (i) can be proved by using the argument in \cite[Proposition 4.6]{AK20}. 

To prove (ii), suppose $\varphi(x,y)$ has $k$-DCTP$_2$ with $(a_\eta)_{\eta\in\tree}$. By Remark \ref{rmk:str-indisc k-DCTP_2}, we may assume that $(a_\eta)_{\eta\in\tree}$ is strongly indiscernible. For each $n<\omega$ and $i<2$, let $b_{n,i}=a_{\lor^n\coc\la 1-i\ra}$. Then $(b_{n,i})_{n<\omega, i<2}$ satisfies the conditions in \cite[Lemma 2.3]{KR20}. Thus $T$ has SOP$_1$.
\end{proof}
\end{rmk}

Therefore, the class of N-$\omega$-DCTP$_2$ theories is a subclass of the class of NATP theories and it is a common extension of the class of NTP$_2$ theories and the class of NSOP$_1$ theories. So we have one more dividing line in the following diagram.

\begin{center}
\includegraphics[width=0.6\linewidth]{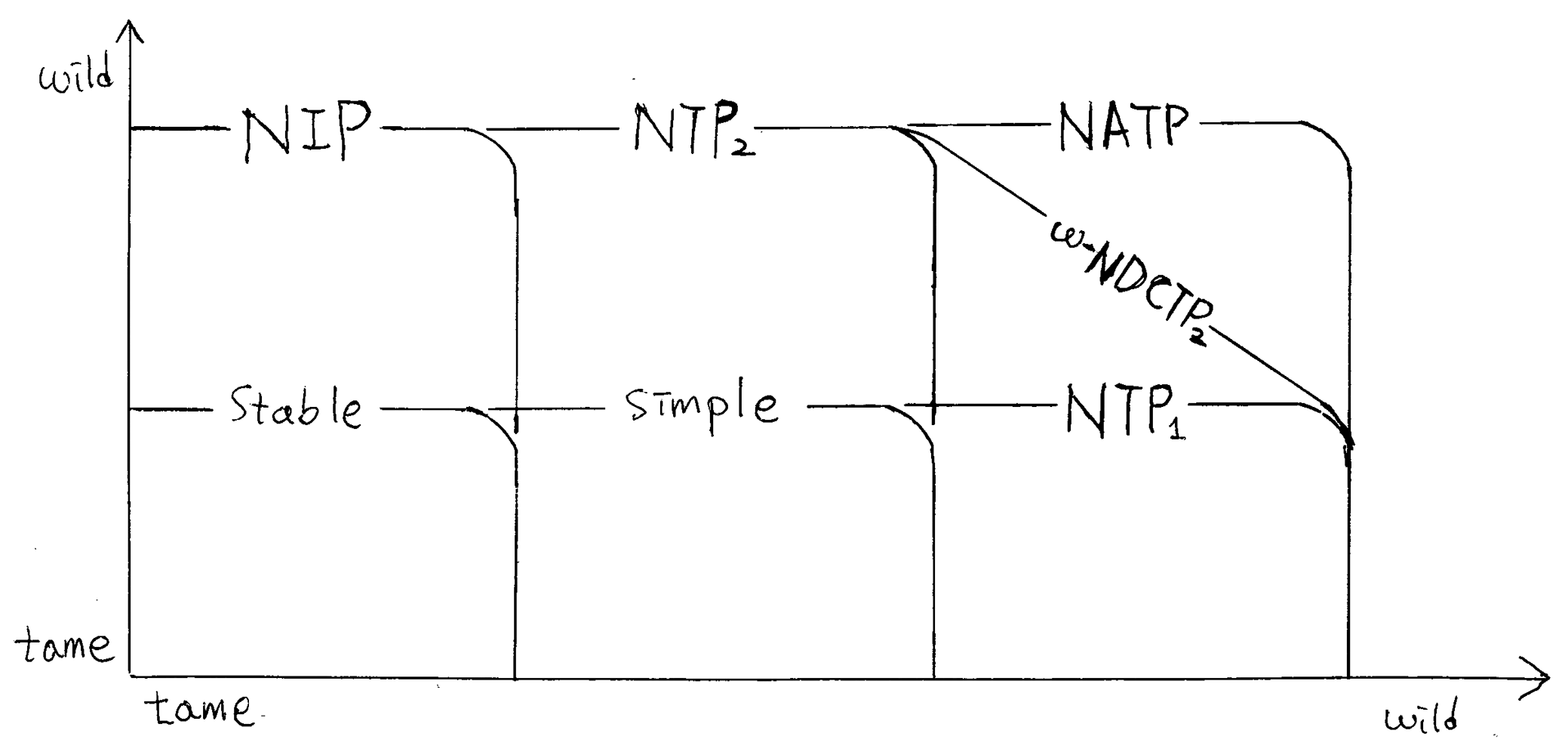}
\end{center}


And in N-$\omega$-DCTP$_2$ theories, it is proved that coheir-dividing and coheir-forking are equivalent over models, as below.

\begin{dfn}\cite{Mut22}\label{def:canonical coheir morley sequences}
 A sequence $(a_i)_{i<\kappa}$ is called a {\it canonical coheir Morley sequence in $\tp(a/M)$} (or in short, we call it a {\it canonical Moerly sequence}) if there exists a global coheir $p(x)\supseteq\tp(a/M)$ over $M$ such that $a_i\models p(x)|_{Ma_{<i}}$ for all $i<\kappa$. We say a formula $\varphi(x,b)$ {\it canonical coheir-divide over $M$} if there exists a canonical coheir Morley sequence $(b_i)_{i<\omega}$ in $\tp(b/M)$ such that $\{\varphi(x,b_i)\}_{i<\omega}$ is inconsistent.
\end{dfn}

\begin{fact}\cite[Theorem 4.9]{Mut22}\label{fact:in k-DCTP2, coheir-dividing=forking}
Suppose $T$ is N-$\omega$-DCTP$_2$. If $\varphi(x,a)$ coheir-divides over a model $M$, then for any canonical Morley sequence $(a_i)_{i<\omega}$ over $M$ with $a_0=a$, 
 the set $\{\varphi(x,a_i)\}_{i<\omega}$ is inconsistent. As a consequence, coheir-dividing and coheir-forking are equivalent over models.
\end{fact}


Using Theorem \ref{thm:Kim-dividing->coheir-dividing}, we can make this result stronger, at least synthetically.

\begin{rmk}
If a theory is N-$\omega$-DCTP$_2$, then Kim-dividing and Kim-forking are equivalent over models.
\begin{proof}
Suppose $\varphi(x,b)$ Kim-forks over $M$, a model. Since the theory is N-$\omega$-DCTP$_2$, it is NATP. By Theorem \ref{thm:Kim-dividing->coheir-dividing}, $\varphi(x,b)$ coheir-forks over $M$. By Fact \ref{fact:in k-DCTP2, coheir-dividing=forking}, $\varphi(x,b)$ coheir-divides over $M$. Thus it Kim-divides over $M$.
\end{proof}
\end{rmk}

\begin{question}
Is there a {\it natural}
example of N-$\omega$-DCTP$_2$ theories having TP$_1$ and TP$_2$?
\end{question}

In the authors' previous work \cite{AKLL23} with JinHoo Ahn and Junguk Lee, we observed that some model theoretic constructions of the first-order theories give examples of NATP theories having TP$_1$ and TP$_2$. It will be interesting if we check whether the similar results hold with N-$\omega$-DCTP$_2$.

\begin{question}
For theories, is having NATP equivalent to having N-$\omega$-DCTP$_2$?
\end{question}


We end this discussion with mentioning relations between Conant-independence and Kim's lemma.

\begin{dfn}\cite{Mut22}
Let $\varphi(x,y)$ be an $\L$-formula, $M$ a model.
\begin{itemize}
\item[(i)] We say {\it $\varphi(x,b)$ Kim-Conant-divides over $M$} if for all invariant Morley sequence $(b_i)_{i<\omega}$ over $M$ starting with $b$, the set $\{\varphi(x,b_i)\}_{i<\omega}$ is inconsistent.
\item[(ii)] We say {\it $\varphi(x,b)$ coheir-Conant-divides over $M$} if for all coheir Morley sequence $(b_i)_{i<\omega}$ over $M$ starting with $b$, the set $\{\varphi(x,b_i)\}_{i<\omega}$ is inconsistent.
\item[(iii)] We say {\it $\varphi(x,b)$ canonical coheir-Conant-divides over $M$ {\rm(in short, we say it} canonical-Conant-divides over $M${\rm)}} if for all canonical coheir Morley sequence $(b_i)_{i<\omega}$ over $M$ starting with $b$, the set $\{\varphi(x,b_i)\}_{i<\omega}$ is inconsistent.
\item[(iv)] We say {\it $a$ and $b$ are Kim-Conant-independent over $M$} and write ${a\ind_M^{Kd^*}}b$ if $\tp(a/Mb)$ has no formula Kim-Conant-dividing over $M$.
\item[(v)] We say {\it $a$ and $b$ are coheir-Conant-independent over $M$} and write ${a\ind_M^{cd^*}}b$ if $\tp(a/Mb)$ has no formula coheir-Conant-dividing over $M$.
\item[(vi)] We say {\it $a$ and $b$ are canonical coheir-Conant-independent over $M$ {\rm (in short,} canonical-Conant-independent over $M${\rm )}} and write ${a\ind_M^{ccd^*}}b$ if $\tp(a/Mb)$ has no formula canonical-Conant-dividing over $M$.
\end{itemize}
\end{dfn}

\begin{rmk}
In N-$\omega$-DCTP$_2$ theories, $\ind^{ccd^*}=\ind^{Kd}=\ind^K$ over models. In fact, the following are equivalent.
\begin{itemize}
\item[(i)] For any model $M$ and $a, b$, if $a\ind^{ccd^*}_Mb$, then $a\ind^{Kd}_Mb$.
\item[(ii)] For any model $M$ and $a, b$, if there exists a coheir Morley sequence $(b_i)_{i<\omega}$ with $b_0=b$ such that $(b'_i)_{i<\omega}$ is not $Ma$-indiscernible for any $(b'_i)_{i<\omega}\equiv_{Mb}(b_i)_{i<\omega}$, then there is no $Ma$-indiscernible canonical Morley sequence $(b_i)_{i<\omega}$ with $b_0=b$.
\item[(iii)] For any model $M$ and $b$, if $\varphi(x,b)$ is consistent and coheir-divides over a model $M$, then the set $\{\varphi(x,b_i)\}_{i<\omega}$ is inconsistent for any canonical Morley sequence $(b_i)_{i<\omega}$ over $M$ with $b_0=b$. 
\end{itemize} 
And (iii) is from Fact \ref{fact:in k-DCTP2, coheir-dividing=forking}.

But as in the example of DLO we mentioned in the introduction, $\ind^{cd^*}$ and $\ind^{K}$ are not equivalent in N-$\omega$-DCTP$_2$ theories in general, even over models.
\end{rmk}

%

\subsection*{Acknowledgement}

We thank Nicholas Ramsey for his helpful comments on Kim-dividing in NSOP$_1$ theories, and Artem Chernikov for letting us know the counterexample in the introduction. We are grateful to Scott Mutchnik for discussing DCTP$_2$ and Bonghun Lee for discussing the subject of this paper in a seminar. We would also like to thank the anonymous referee for the valuable comments and suggestions.

The authors were supported by NRF of Korea grants 2021R1A2C1009639. The first author is supported by a KIAS individual grant, project no. 6G091801.
The second author is supported by the 2023, 2024 Yonsei University post-doc. researcher supporting program, project no. 2023-12-0159, 2024-12-0214, and the G-LAMP Program of the NRF of Korea grant funded by the Ministry of Education(No. RS-2024-00441954).

\end{document}